\documentclass[preprint,12pt]{elsarticle}

\usepackage{amssymb}
\usepackage{amsthm}
\usepackage{amsmath}
\usepackage{siunitx}
\usepackage{tabularx}
\usepackage{bm}
\usepackage{color}
\usepackage{moreverb}
\graphicspath{ {./Figures/} }
\usepackage[utf8]{inputenc}
\usepackage{psfrag}
\theoremstyle{plain}

\theoremstyle{plain}

\theoremstyle{definition}
\newtheorem{remark}{Remark}
\theoremstyle{definition}

\theoremstyle{definition}

\usepackage{amsfonts}
\usepackage{float}
\usepackage{multirow}
\usepackage{setspace}
\usepackage{flushend}
\usepackage{array}
\usepackage{tabu}
\usepackage[dvipsnames]{xcolor}
\usepackage[english]{babel}
\usepackage{lipsum} 
\usepackage{url}
\usepackage{dsfont}
\usepackage[figuresright]{rotating}
\usepackage{subfigure}
\usepackage{svgcolor}
\usepackage{makecell}
\usepackage{booktabs}
\usepackage{scalerel}
\usepackage{eurosym}
\usepackage{footnote}
\usepackage{cancel}
\usepackage{caption}
\usepackage[dvips]{epsfig}
\usepackage{graphicx}

\usepackage{subcaption}

\newtheorem{prop}{Proposition}
\newtheorem{proper}{Property}

\journal{Ocean Engineering}

\begin{document}

\begin{frontmatter}

%% Title, authors and addresses

%% use the tnoteref command within \title for footnotes;
%% use the tnotetext command for theassociated footnote;
%% use the fnref command within \author or \address for footnotes;
%% use the fntext command for theassociated footnote;
%% use the corref command within \author for corresponding author footnotes;
%% use the cortext command for theassociated footnote;
%% use the ead command for the email address,
%% and the form \ead[url] for the home page:
%% \title{Title\tnoteref{label1}}
%% \tnotetext[label1]{}
%% \author{Name\corref{cor1}\fnref{label2}}
%% \ead{email address}
%% \ead[url]{home page}
%% \fntext[label2]{}
%% \cortext[cor1]{}
%% \affiliation{organization={},
%%             addressline={},
%%             city={},
%%             postcode={},
%%             state={},
%%             country={}}
%% \fntext[label3]{}

\title{Practical identification approach for the actuation dynamics of autonomous surface vehicles with minimal instrumentation: extended version}

%% use optional labels to link authors explicitly to addresses:
%% \author[label1,label2]{}
%% \affiliation[label1]{organization={},
%%             addressline={},
%%             city={},
%%             postcode={},
%%             state={},
%%             country={}}
%%
%% \affiliation[label2]{organization={},
%%             addressline={},
%%             city={},
%%             postcode={},
%%             state={},
%%             country={}}

\author[inst2]{Thalia Morel}
\author[inst1,inst3]{Luis Orihuela}
\author[inst4]{Christophe Combastel}
\author[inst2]{Guillermo Bejarano}

\affiliation[inst2]{organization={Department of Engineering, Universidad Loyola Andalucía},%Department and Organization
            addressline={Avenida de las Universidades, 2}, 
            postcode={41704}, 
            city={Dos Hermanas (Seville)},
            country={Spain}}

\affiliation[inst1]{organization={Dept. Ingeniería Electrónica, Sistemas Informáticos y Automática, Universidad de Huelva},%Department and Organization
            addressline={Avda. de las Fuerzas Armadas, no number}, 
            postcode={21007}, 
            city={Huelva},
            country={Spain}}

\affiliation[inst3]{organization={Centro de Investigación en Tecnología, Energía y Sostenibilidad, Universidad de Huelva},%Department and Organization
            addressline={Ctra. Huelva-Palos de la Frontera, Campus La Rábida}, 
            postcode={21819}, 
            city={ Palos de la Frontera (Huelva)},
            country={Spain}}            

\affiliation[inst4]{organization={Univ. Bordeaux, CNRS, Bordeaux INP, IMS, UMR 5218},%Department and Organization
            %addressline={}, 
            postcode={F-33400}, 
            city={Talence},
            country={France}}

\begin{abstract}
A practical method for identifying the propeller model and inertia matrix of a marine Autonomous Surface Vehicle (ASV) is proposed in this work. Special attention is paid to limiting the instrumentation requirements. Based on a generic grey-box dynamic modelling addressing the considered catamaran-shaped ASV architecture, the static/dynamic behaviour of both propellers and the vessel dynamic are jointly estimated using the sole measurements of position, heading, and propellers pulse width modulation (PWM) signals. No accelerometer is required. Two distinct grey-box configurations involving either a static polynomial or a dynamic modelling of each propeller are proposed and compared. The resulting ASV identification methodology is shown to provide insight into the whole vessel inertial characteristics, which are key enablers in the development of autonomous navigation and control systems. Model validation was performed using data collected from the reported experiments. Model prediction errors related to both linear velocities and yaw rate are evaluated and compared based on given metrics. The results underscore the robustness and accuracy of the identified models in capturing the essential dynamics of the ASV, with a determination coefficient that consistently exceeds 0.94 for all estimated velocities.
\end{abstract}

% %%Graphical abstract
% \begin{graphicalabstract}
% \includegraphics{grabs}
% \end{graphicalabstract}

% %%Research highlights
% \begin{highlights}
% \item Research highlight 1
% \item Research highlight 2
% \end{highlights}

\begin{keyword}
%% keywords here, in the form: keyword \sep keyword
Autonomous surface vehicles \sep System identification \sep Actuation modelling \sep Inertia \sep Vessel dynamics

%System identification \sep Actuation modelling \sep Autonomous surface vehicles
% %% PACS codes here, in the form: \PACS code \sep code
% \PACS 0000 \sep 1111
%% MSC codes here, in the form: \MSC code \sep code
%% or \MSC[2008] code \sep code (2000 is the default)
\MSC[2020] 93B30 \sep 93C10
\end{keyword}

\end{frontmatter}

%% \linenumbers

%% main text
%%%%%%%%%%%%%%%%%%%%%%%%%%%%%%%%%%%%%%%%%%%%%%%%%
\section{Introduction}
\label{sec:Introduction}
%%%%%%%%%%%%%%%%%%%%%%%%%%%%%%%%%%%%%%%%%%%%%%%%%

Recent advances in ASVs have opened up new and exciting possibilities in various applications, offering a potential reduction in the reliance on large manned ships that cause high operational costs~\cite{jovanovic2024review, gantiva2024coordinated}. These autonomous vehicles have the ability to revolutionise several fields, including search and rescue operations, ocean exploration/monitoring, border surveillance, maritime security among other civilian and military applications. Using ASVs, these tasks can be performed efficiently, leading to improved safety, cost savings, and improved environmental stewardship~\cite{morel2022experimental, gantiva2024implementation}.

ASVs have gained significant attention for both research and civil applications since the early 1990s~\cite{simetti2022control}. Notable prototypes include the catamaran~\emph{AutoCat} for the collection of hydrographic data~\cite{wu2022identification}, the autonomous kayak \emph{SCOUT} for acoustic navigation and autonomous system tests~\cite{curcio2005scout}, and the \emph{Measuring Dolphin} of the University of Rostock for water monitoring~\cite{MeasuringDolphin2006}. Other examples include the \emph{Charlie} catamaran from CNR-ISSIA for sea surface sampling~\cite{caccia2007unmanned} and the \emph{Delfim} from IST-ISR for enhancing acoustic communication~\cite{alves2006vehicle}.

With a similar purpose of water body monitoring and pollution tracking, an autonomous catamaran-like surface vessel named \emph{Yellowfish} has been conceived and developed at Universidad Loyola Andalucía. This vessel is designed to navigate through calm water and at low speed, due to its size and the power of its actuation system, while relying on a minimal instrumentation scheme.

Motion control is crucial to ensure safe navigation and successful path following of the ASV, despite the significant effect of non-measured disturbances due to wind, waves, and currents in real environments~\cite{breivik2010topics}. Most control design strategies presented in the literature rely on an open-loop dynamic model of the vessel behaviour including its actuation system. However, obtaining an accurate model for an ASV based on physical principles can be challenging due to strong non-linearities, unmodelled hydrodynamics related to neglected degrees of freedom, environmental disturbances, and internal parametric uncertainty. However, modelling and identifying such systems remain essential for accurate navigation and control of the vessel~\cite{peng2020overview}. Conventional approaches often require extensive instrumentation, including acceleration measurements, which may increase both the cost and complexity of the system.

Several experimental works in the literature have addressed the identification of dynamic models for the motion of different kinds of ASVs. In the literature, the works analysed in the following differ in their approach, objectives and layout of the vehicles considered, and they give insight into the current state of the art. This contributes to highlight the main techniques used, as well as some of their main restrictions to \emph{jointly} identify both the ASV dynamics and actuation system.

Sonnenburg and Woolsey~\cite{Eriksen2017} developed a planar motion model for a single-hull ASV, comparing experimentally identified models under different speeds and conditions, using a modified rigid-hull inflatable boat~\cite{Sonnenburg2013}. Their multiple-model approach adjusts based on nominal speed. This paper proposes a control-oriented modelling procedure for a high-speed rigid-hull inflatable ASV across various speed regions, without explicitly modelling the actuation system.

Wu~\emph{et al.} proposed a method to identify linear and non-linear hydrodynamic coefficients using a Support Vector Machine (SVM)~\cite{wu2022identification} for the modelling of ASV manoeuvres. Linear coefficients were identified from turning tests at small water-jet angles, while larger angles and the linear coefficients were used for non-linear coefficient identification. The model accurately describes the manoeuvrability of the ASV. Xu~\emph{et al.}. combined the Least Square Support Vector Machine (LSSVM) with the Cuckoo Search (CS) algorithm to identify a 3-DOF dynamic model of an ASV~\cite{xu2020identification}, outperforming other methods to predict surge and sway velocities, although accurate velocity measurements are required.

In~\cite{abrougui2021modeling}, a 3-DOF control-oriented dynamic model for an ASV was proposed, with parameters determined through theoretical and experimental methods. The control system enables the ASV to follow a desired trajectory accurately, leading to the conclusion that lateral sway is negligible in closed-loop control. The emphasis is on achieving sufficient accuracy for motion control rather than precise parameter estimation.

Caccia \emph{et al.} proposed a least-squares identification procedure for ASVs using basic sensors onboard that are typically found in small, low-cost vessels~\cite{caccia2008practical}. Validation was performed using data from auto-heading and line-following manoeuvres in typical conditions. Linear and quadratic models were applied to the propellers in the actuation system. The linear drag model performed better at low angular rates but overestimated yaw rates at high rudder angles, unlike the quadratic model.

Simetti and Indiveri recently proposed a first-principle manoeuvring model for a small twin-thruster ASV, addressing non-linearity and asymmetry due to the catamaran layout's unequal port/starboard dynamics \cite{simetti2022control}. The model adapts hydrodynamic derivatives for large or tight turns and is identified by correlating constant thruster revolutions per minute (RPM) with steady-state surge and yaw rates from the Global Navigation Satellite System (GNSS) and inertial data. They propose identifying both the ASV's dynamic and actuation models concurrently, but sway dynamics are not identified, and the damping matrix is assumed to be diagonal.

As a summary of the literature review, despite the different approaches, focus, and ASV layouts, three main facts appear:
\begin{itemize}
    \item Previous works address either the identification of the dynamic model of the vessel, or the actuation model, but not the combined effect on the vehicle motion, at least without strong assumptions on the model structure or the availability of accurate velocity measurements.   
    \item Those works that focus on the identification of simplified control-oriented dynamic models sacrifice, in most cases, the identification of individual parameters of the model.
    \item Concerning the description of the actuation system, both linear and quadratic static models have been considered. Dynamic effects are omitted in most cases.
\end{itemize}

%% Description of the contributions:  
In this work, a practical identification approach is proposed to obtain a
control-oriented dynamic model of an ASV with limited instrumentation. Reported \emph{Yellowfish} operations provided the experimental data sets (made available) used to validate the proposed methodology. It focuses on the identification of the so-called \emph{input gain} term~\cite{peng2020output}, which relates the PWM signals applied to the propellers to the accelerations generated in the vessel, thus merging the actuation system with the inertial properties of the vessel. This term turns out to be key for motion control and state estimation~\cite{ abrougui2021modeling, peng2020output}. In this context, the main contribution of the paper consists in proposing an original grey-box modelling approach involving polynomial and dynamic terms to accurately identify this so-called \emph{input gain}, based on a minimal instrumentation scheme. The comparison and experimental validation of the resulting models are shown to provide valuable insights not only into the characteristics of both the propellers and the inertia matrix but also into the whole ASV dynamic.

%% Document organization:
The paper is organised as follows. After this introduction, the problem is stated in Section~\ref{sec:sample:ProblemStatement}. Section~\ref{sec:sample:DescriptionASV} describes the experimental platform, detailing the \emph{Yellowfish} ASV with the instrumentation used and the available data. Section~\ref{sec:sample:Modelling} details the models considered for both the vessel dynamics and the actuation system. Then, Section~\ref{sec:sample:Identification} describes the proposed identification procedure for all actuation models considered. Section~\ref{sec:sample:Results} describes the experimental results, presents the obtained models, and reports the validation results. It also discusses model selection, considering the uncertainty in identifying model coefficients. Finally, concluding remarks are given in Section~\ref{sec:sample:Conclusions}.

%%%%%%%%%%%%%%%%%%%%%%%%%%%%%%%%%%%%%%%%%%%%%%%%%
\section{Problem Statement} 
\label{sec:sample:ProblemStatement}
\begin{figure}[h]
    \begin{centering}
    \includegraphics[angle=0,width=0.5\textwidth]{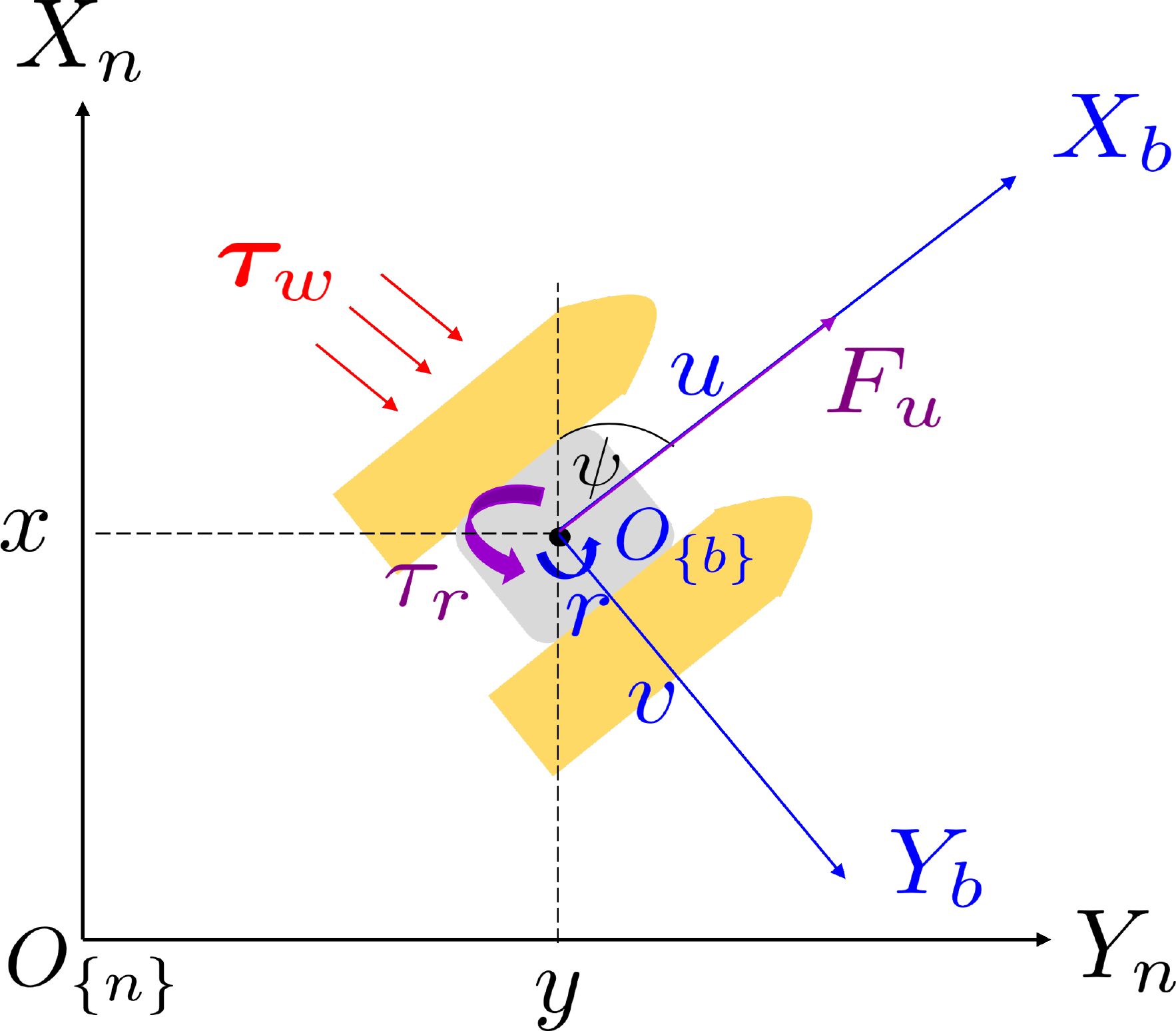}
    \par\end{centering}
    \caption{Reference Frames of the Vessel}
    \label{fig_ASVcoord}
\end{figure}

%%%%%%%%%%%%%%%%%%%%%%%%%%%%%%%%%%%%%%%%%%%%%%%%%
In order to state the ASV grey-box modelling and identification problem addressed in this work, the kinematics and kinetics of a surface vessel in horizontal motion are first described using the following dynamic model~\cite{fossen2011handbook}:

\begin{equation}
	\begin{aligned}
		\left\{ 
		\begin{matrix}
			\begin{aligned}
				\bm{\dot{\eta}} &= \bm{R}(\psi) \bm{\nu} \\
				\bm{M} \bm{\dot{\nu}} &= - \bm{C}(\bm{\nu})\bm{\nu} - \bm{D}(\bm{\nu})\bm{\nu} - \bm{g}(\bm{\eta}) + \bm{\tau_{w}} + \bm{\tau} , \\
			\end{aligned}
		\end{matrix}
		\right.
	\end{aligned}
	\label{eqFossenModel}
\end{equation}
where~$\bm{\eta} = [x \; y \; \psi]^{T}$ is the planar position and heading vector expressed in the earth-fixed inertial frame \{n\} (see Figure~\ref{fig_ASVcoord}), $\bm{\nu} = [u \; v \; r]^{T}$ is the linear and angular velocity vector expressed in the body-fixed frame \{b\}, $\bm{\tau}~=~[F_u \; F_v \; \tau_r]^{T}$ refers to the force/torque vector, which represents the control action, and~$\bm{\tau}_{w} = \left[ F_{w,u} \; F_{w,v} \; \tau_{w,r} \right]^{T}$ refers to the environmental forces and torque due to wind, waves, and currents. $\bm{M}$ is the inertia matrix, $\bm{C}(\bm{\nu})$ the Coriolis and centrifugal matrix, $\bm{D}(\bm{\nu})$ the damping matrix, and~$\bm{g}(\bm{\eta})$ describes the gravitational and buoyancy forces, which are neglected for a surface vehicle in planar movements. Note that an underactuated actuation system is considered, i.e. $\bm{\tau}~=~[F_u \; 0 \; \tau_r]^{T}$, since there is no lateral force in $\bm{\tau}$. 

As first stated in~\cite{liu2019state}, the lumped generalised disturbance vector~$\bm{\sigma} = \left[ \sigma_u \; \sigma_v \; \sigma_r \right]^T\in\mathbb{R}^3$ can be defined as:
\begin{equation}
	\begin{aligned}
		\bm{\sigma} \equiv \bm{M}^{-1}\left( -\bm{C}(\bm{\nu})\bm{\nu} - \bm{D}(\bm{\nu})\bm{\nu} - \bm{g}(\bm{\eta}) + \bm{\tau_w} \right) . \\
	\end{aligned} 
	\label{eqLiuSigma}
\end{equation}

Vector~$\bm{\sigma}$ groups together disturbances due to hydrodynamic effects, internal unmodelled and coupling dynamics, uncertainties, and external disturbances. Using~\eqref{eqLiuSigma}, the body dynamic model~\eqref{eqFossenModel} can be redefined as:
\begin{equation}
	\begin{aligned}
		\left\{ 
		\begin{matrix}
			\begin{aligned}
				\bm{\dot{\eta}} &= \bm{R}(\bm{\psi}) \bm{\nu}, \\
				\bm{\dot{\nu}} &= \bm{M}^{-1}\bm{\tau} + \bm{\sigma}. 
			\end{aligned}
		\end{matrix}
		\right.
	\end{aligned}
	\label{eqFossenModel_Liu}
\end{equation} 

The so-called input gain $\bm{G}~=~[G_u \; G_v \; G_r]^{T} \in \mathbb{R}^{3}$ is then defined as:
\begin{equation}
    \bm{G} \equiv \bm{M}^{-1}\bm{\tau},
    \label{eqkineticIG}
\end{equation}

\noindent and depends on both the vessel inertial characteristics and the truly applied propelling force and torque. Though the latter are closely related to the propellers' PWM control signals (denoted as $\delta_L$ and $\delta_R$ for the left and right propeller, respectively), the corresponding relation $\bm{G}(\delta_L, \delta_R)$ is most often far from trivial. In fact, it may feature nonlinear and dynamic effects, including possible switching related to propellers' forward and reverse modes. Moreover, a simple characterisation from first physical principles (knowledge model) cannot be easily derived from the sole propeller characteristics since the truly applied propelling force and torque depend not only on the position of propellers on the vessel but also on much more complex and possibly coupled hydrodynamics phenomena. This motivates the development of an identification methodology primarily focused on the estimation of $\bm{G}(\delta_L, \delta_R)$ from experimental data. 

More precisely, the aim of this work is,
\begin{itemize}
    \item firstly, to define and compare relevant grey-box models (non-linear, dynamic) for the so-called input gain $\bm{G}=\bm{G}(\delta_L, \delta_R)$ in~\eqref{eqkineticIG}, 
    \item then estimate all related parameters and, 
    \item  finally, validate the resulting whole ASV dynamic model based on explicitly given metrics.
\end{itemize}
Another requirement of the identification methodology is to rely on a minimal instrumentation scheme reduced to the linear position measurements given by the GNSS system, the heading measured by the Attitude and Heading Reference Systems (AHRS), and the propellers' PWM control signals. Such signals are easily made available for a wide class of ASVs, including the catamaran-shaped \emph{YellowFish} described in Section~\ref{sec:sample:DescriptionASV} and serving as an illustrative application case in this paper.

%%%%%%%%%%%%%%%%%%%%%%%%%%%%%%%%%%%%%%%%%%%%%%%%%
\section{Experimental Platform, Instrumentation and Available Data}
\label{sec:sample:DescriptionASV}
%%%%%%%%%%%%%%%%%%%%%%%%%%%%%%%%%%%%%%%%%%%%%%%%%

\begin{figure}[h]
    \begin{centering}
    \includegraphics[angle=0,width=0.8\textwidth]{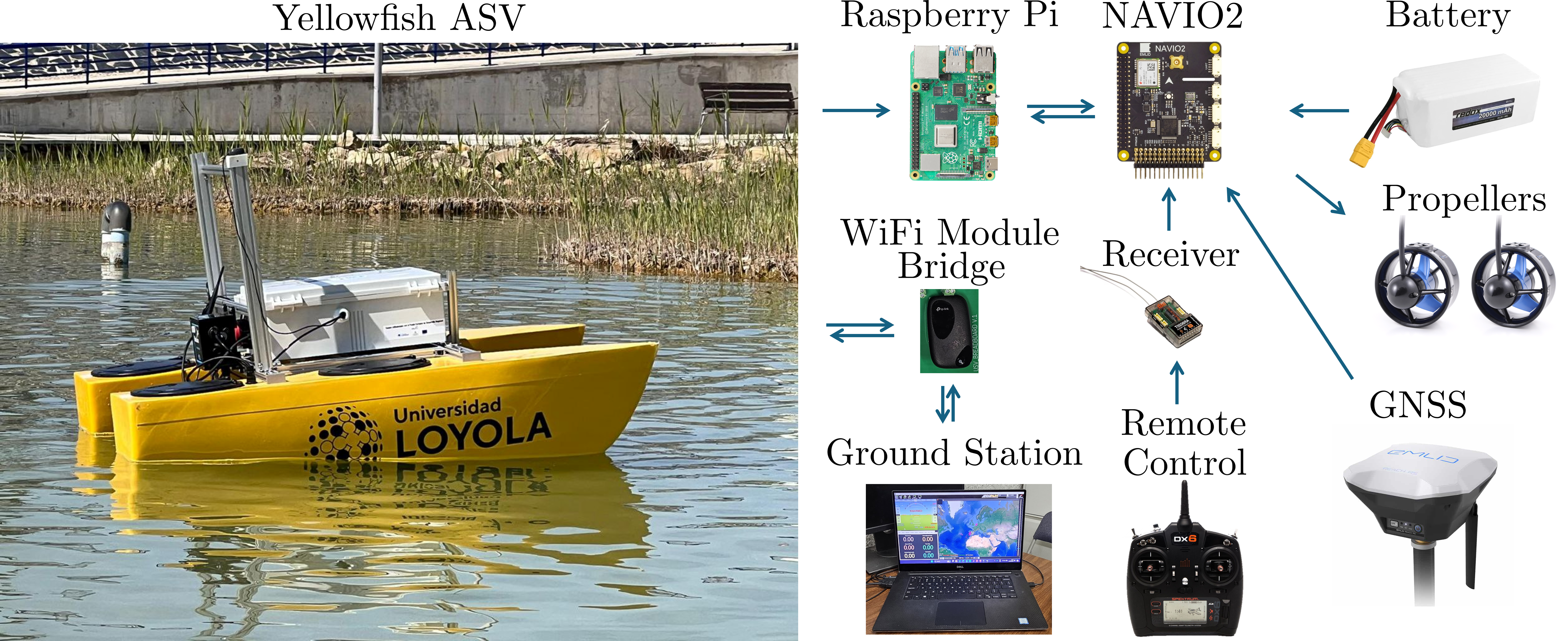}
    \par\end{centering}
    \caption{The \emph{Yellowfish} ASV Experimental Platform}
    \label{fig_ASVplatform}
\end{figure}

The \emph{Yellowfish} ASV, depicted in Figure~\ref{fig_ASVplatform}, is a catamaran vessel made of fibreglass, measuring $\SI{128}{\centi\meter}$ long and $\SI{98}{\centi\meter}$ wide. The catamaran design features two identical hulls, firmly connected by aluminium bars. 
The propulsion of the \emph{Yellowfish} ASV is provided by two T200 thrusters manufactured by BlueRobotics, which are electric propeller motors.

Regarding hardware, the ASV features a Raspberry Pi 4 model B integrated with a Navio2 for command and data processing. The Raspberry Pi 4, a mini-computer, operates with an integrated system specifically tailored to align with the characteristics of Navio2. Its primary role within the system is to facilitate telemetry support. The experimental platform, equipped with the necessary sensors and hardware for data acquisition, is shown in Figure~\ref{fig_ASVplatform}.

Furthermore, the software implemented on the Navio2 is called Ardupilot. Essentially, Navio2 uses Ardupilot for comprehensive sensor configuration and vehicle control across all installed vehicles. Mission Planner serves as the ground station software utilised for telecommunication purposes with the ship. The DataFlash logs can be obtained through Mavlink from the Mission Planner and create a MATLAB file with all sensor values. For additional details on the \emph{Yellowfish} ASV's hardware and software, the reader is referred to~\cite{morel2022modelling,morel2023modelling}.

The \emph{Yellowfish} ASV mounts different sensors whose available data and their respective sampling times are presented in Table~\ref{tab_ASV_variables}. As seen in the table, the sensors have different sampling times: The AHRS at $20$ ms and the GNSS at $200$ms. In addition to these outputs, the PWM signals generated with the radio controller (RCOU) can be measured with a sampling time of $100$ms. 

Based on the data provided by the mentioned sensors, some transformations are necessary to obtain the variables used in the vessel model. First, each propeller (left and right) takes a PWM signal as input. The PWM values of the left and right thrusters are normalised between $[-1, 1]$. This signal, after normalisation, can be described by $\delta_i$ $\in$ $[-1, 1]$, $i\in\{L,R\}$, being $\delta_i=\pm1$ a PWM signal with a duty cycle of $100\%$, and $\delta_i=0$ a PWM signal with a duty cycle of $0\%$. More details on the normalisation of the PWM signals can be found in~\cite{morel2022modelling}. The propellers produce forward motion forces for positive values of $\delta_i$. Other metrics, such as planar positioning coordinates \(\{x, y\}\), are obtained by transforming point locations provided by the GNSS from geographic coordinates to local Cartesian coordinates, returned as north-east-down (NED). After this transformation, a location shift is performed, as the GNSS antenna is not positioned at the same place as the Navio2, from which other data are derived.  The heading $\psi$ is measured directly by the AHRS.

\begin{table}[h]
\centering
\caption{Summary of Sensor Data and Variables for Vessel Model Analysis}
\label{tab_ASV_variables}
\resizebox{\columnwidth}{!}{%
\begin{tabular}{llcccl}
\toprule
\textbf{Sensor} & \textbf{Available Data} & \textbf{Sampling Rate} & \textbf{Obtained Variable} & \textbf{Unit} & \textbf{Description of Variable} \\ \hline
GNSS            & Latitude, longitude     & 200ms                  & $x,y$                        & m             & North and east absolute position \\
RCOU            & $\text{PWM}_L$, $\text{PWM}_R$             & 100ms                  & $\delta_L, \delta_R$            & $\mu$s            & Propellers normalised duty cycle            \\
AHRS            & Heading                & 20ms                   &$\psi$                         & rad           & Heading   angle                       \\ \toprule
\end{tabular}%
}
\end{table}

%%%%%%%%%%%%%%%%%%%%%%%%%%%%%%%%%%%%%%%%%%%%%%%%%
\section{Modelling} 
\label{sec:sample:Modelling}
\begin{figure}[ht]
    \begin{centering}
    \includegraphics[angle=0,width=\textwidth]{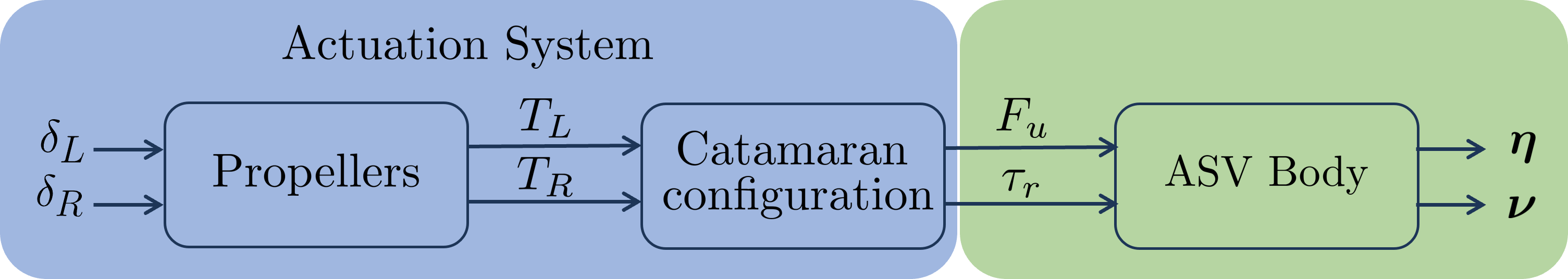}
    \par\end{centering}
    \caption{Block Diagram of the ASV, Highlighting the Actuation System and the ASV Body}
    \label{fig_ASVmodel}
\end{figure}

%%%%%%%%%%%%%%%%%%%%%%%%%%%%%%%%%%%%%%%%%%%%%%%%%
This section aims to find the expression for the input gain as a grey-box model, requiring the complete modelling of the vehicle. The overall vehicle model is represented by the block diagram in Figure~\ref{fig_ASVmodel}, where the system is divided into two primary components: the actuation system and the ASV body.
%The modelling of the complete vehicle can be described with the block diagram shown in Figure~\ref{fig_ASVmodel}. The system is divided into two main components: the actuation system and the ASV body. 
The actuation system takes the PWM signals for the left and right propellers, denoted as $\delta_L$ and $\delta_R$, as inputs. It then produces the force $F_u$ and the torque $\tau_r$ applied to the vessel body as outputs, using thrusts $T_L$ and $T_R$ as intermediate variables within the actuation system. The ASV body takes the force $F_u$ and torque $\tau_r$ as inputs and produces, as outputs, the full state of the vessel, consisting of its 3-DOF position $\bm{\eta}$ and velocity $\bm{\nu}$ vectors. In the following subsections, further details are given for the models of each block.

\subsection{ASV body}

Based on the kinetic equation presented in~\eqref{eqFossenModel_Liu}, considering manoeuvring theory (almost constant positive forward speed and relatively calm water) and assuming homogeneous mass distribution and~$X_b-Z_b$-plane symmetry, being $X_b$ and $Z_b$ the aft-to-fore and top-to-bottom axes, respectively (see~\cite{fossen2011handbook}), the structure of~$\bm{M}$ is shown in~\eqref{eqMatrixM}:
\begin{equation}
   \begin{aligned}
	\bm{M} &= \left[
	\begin{matrix}
		m - X_{\dot{u}} 		& 0 					& 0 	 				\\
		0 						& m - Y_{\dot{v}} 		& m x_g - Y_{\dot{r}} 	\\
		0 						& m x_g - Y_{\dot{r}}   & I_z - N_{\dot{r}}  	\\
	\end{matrix}
	\right],
	\label{eqMatrixM}
	\end{aligned}
\end{equation}
where~$m$ is the ASV total mass, $x_g$ is the distance from the centre of gravity of the vessel to the origin of \{b\} (measured along~$X_b$), 
$I_z$ is the moment of inertia about the~$Z_{b}$ axis, and~$X_{(\cdot)}$, $Y_{(\cdot)}$, and $N_{(\cdot)}$ are hydrodynamic parameters according to standard notation~\cite{sname1950nomenclature}. Notice that $\bm{M} = \bm{M}^T$, and the structure of $\bm{M}$ ensures that $\bm{M}^{-1}$ exists for well-designed vessels \cite{fossen2011handbook}.

According to~\cite{fossen2011handbook}, Coriolis $\bm{C}(\bm{\nu})$ and damping $\bm{D}(\bm{\nu})$ matrices are given by:
\scriptsize
\begin{align}
\bm{C}(\bm{\nu})&=\left[\begin{array}{ccc}
0 & 0 & -(m-Y_{\dot{v}})v-(mx_{g}-Y_{\dot{r}})r\\
0 & 0 & (m-X_{\dot{u}})u \\
(m-Y_{\dot{\nu}})\nu+(mx_{g}-Y_{\dot{r}})r & -(m-X_{\dot{u}})u & 0
\end{array}\right], \label{eqMatrixC}\\
\bm{D}(\bm{\nu})&=\left[\begin{array}{ccc}
-X_{u}-X_{|u|u}|u| & 0 & 0 \\
0 & -Y_{v}-Y_{|v|v}|v|-Y_{|r|v}|r| & -Y_{r}-Y_{|v|r}|v|-Y_{|r|r}|r| \\
0 & -N_{v}-N_{|v|v}|v|-N_{|r|v}|r| & -N_{r}-N_{|v|r}|v|-N_{|r|r}|r|
\end{array}\right].\label{eqMatrixD}
\end{align}
\normalsize

Note that the damping matrix considers non-linear terms, which are neglected in most identification papers but which will be kept in this work.

The rotation matrix~$\bm{R}(\psi)$ between  \{b\} and \{n\} is given by:
\begin{equation} 
	\bm{R}(\psi) = \left[
	\begin{matrix}
		\text{cos}(\psi) & -\text{sin}(\psi) & 0 \\
		\text{sin}(\psi) &  \text{cos}(\psi) & 0 \\
		0		  &	 0		   & 1 \\
	\end{matrix}
	\right] = 
	\left[
	\begin{matrix}
		\bm{R}_2(\psi) & \bm{0} \\
		\bm{0}	 & 1 \\
	\end{matrix}
	\right],
	\label{eqMatrixR}
\end{equation}
where~$\bm{R}_2(\psi) \in \mathbb{R}^{2\times2}$ is the upper left submatrix of~$\bm{R}(\psi)$.

In this study, the identification process relies on a discretized version of the kinetics equation, chosen for the incorporation of unknown disturbances. The discrete-time kinetics equation is then derived, with $h$ denoting the sampling period, as follows:
\begin{equation}
\bm{\nu}(k+1) = \bm{\nu}(k)+ h\bm{M}^{-1}\bm{\tau}(k) + h\bm{\sigma}(k).
\label{eq_discretized_dynamics}
\end{equation}

More precisely, the discretized system dynamics presented in~\eqref{eq_discretized_dynamics} can be divided into three equations, corresponding to every velocity component, as shown: 
\begin{equation}
    \begin{aligned}
        u(k+1) &=u(k)+h\bm{M}^{-1}(1,1)F_u(k) + h\sigma_u(k), \\
        v(k+1) &= v(k)+h\bm{M}^{-1}(2,3)\tau_r(k) + h\sigma_v(k), \\
        r(k+1) &=r(k)+h\bm{M}^{-1}(3,3)\tau_r(k) + h\sigma_r(k), \\
    \end{aligned}
    \label{eq_nu_definition}
\end{equation}
where $\bm{M}^{-1}(i,j)$ stands for the $(i,j)$ element of the inverse of the inertia matrix. Then, it is evident that the input gains can be mathematically defined as:
\begin{equation}
    \begin{aligned}
        G_u(k)&:= h\bm{M}^{-1}(1,1)F_u(k),\\
        G_v(k)&:= h\bm{M}^{-1}(2,3)\tau_r(k), \\
        G_r(k)&:= h\bm{M}^{-1}(3,3)\tau_r(k). \\
    \end{aligned}
    \label{eq_IG_definition2}
\end{equation}

From~\eqref{eqLiuSigma}, it can be assumed that each of the components of $\bm{\sigma}$ can be approximated as a quasi-quadratic function of $\bm{\nu}$ plus some disturbances:
\begin{equation}\label{eqAproxSigma}
    \sigma_{j} = \bm{\nu}^{T}\bm{P_{j}}|\bm{\nu}|+\bm{\nu}^{T}\bm{Q_{j}}\bm{\nu}+\bm{R_{j}}\bm{\nu}+c_{j},
\end{equation}
where $j\in\{u,v,r\}$, $\bm{P}_{j}\in\mathbb{R}^{3\times 3}$, $\bm{R}_{j}\in\mathbb{R}^{1\times 3}$, and $\bm{Q}_{j}\in\mathbb{R}^{3\times 3}$ are constant matrices assumed to be symmetric, and $c_{j}\in\mathbb{R}$ is a possible time-varying vector accounting for the effect of wind, waves, and currents. In particular:
\small
\begin{equation*}
    \bm{P}_u = \left[\!\!\begin{array}{ccc}
    P_u^{(1,1)} & 0 & 0 \\
    0 & 0 & 0 \\
    0 & 0 & 0
    \end{array}\!\!\right],\;\bm{Q}_u = \left[\!\!\begin{array}{ccc}
    0 & 0 & 0 \\
    0 & 0 & Q_{u}^{(2,3)} \\
    0 & Q_{u}^{(2,3)} & Q_{u}^{(3,3)}
    \end{array}\!\!\right],\;\bm{R}_u = \left[\!\!\begin{array}{ccc}
    R_{u}^{(1)} & 0 & 0 
    \end{array}\!\!\right]
\end{equation*}
\begin{equation*}
    \bm{P}_v = \left[\!\!\begin{array}{ccc}
    0 & 0 & 0 \\
    0 & P_v^{(2,2)} & P_v^{(2,3)} \\
    0 & P_v^{(3,2)} & P_v^{(3,3)}
    \end{array}\!\!\right],\;\bm{Q}_v = \left[\!\!\begin{array}{ccc}
    0 & Q_{v}^{(1,2)} & Q_{v}^{(1,3)} \\
    Q_{v}^{(1,2)} & 0 & 0 \\
    Q_{v}^{(1,3)} & 0 & 0
    \end{array}\!\!\right],\;\bm{R}_v = \left[\!\!\begin{array}{ccc}
    0 & R_{v}^{(2)} & R_{v}^{(3)} 
    \end{array}\!\!\right]
\end{equation*}
\begin{equation*}
    \bm{P}_r = \left[\!\!\begin{array}{ccc}
    0 & 0 & 0 \\
    0 & P_r^{(2,2)} & P_r^{(2,3)} \\
    0 & P_r^{(3,2)} & P_r^{(3,3)}
    \end{array}\!\!\right],\;\bm{Q}_r = \left[\!\!\begin{array}{ccc}
    0 & Q_{r}^{(1,2)} & Q_{r}^{(1,3)} \\
    Q_{r}^{(1,2)} & 0 & 0 \\
    Q_{r}^{(1,3)} & 0 & 0
    \end{array}\!\!\right],\;\bm{R}_r = \left[\!\!\begin{array}{ccc}
    0 & R_{r}^{(2)} & R_{r}^{(3)} 
    \end{array}\!\!\right]
\end{equation*}
\normalsize

\subsection{Actuation System}

The physical actuation system comprises two parts. First, the propellers that receive individual normalised PWM signals $\delta_L$, $\delta_R$ and generate thrusts $T_L(\delta_L)$, $T_R(\delta_R)$. This first part has received a lot of attention in the literature, as propellers have complex nonlinear dynamics and are difficult to model and identify~\cite{simetti2022control,blanke2000dynamic,mu2018modeling,zuev2020fault,simmons2021system}. In the current configuration, it is reasonable to assume that both propellers are identical. Thus, a unique function $T(\delta_i),i\in\{L,R\},$ is considered valid for both propellers.

The second part is the particular configuration of the vessel that allows a static transformation from the thrusts of each propeller to the force and torque received by the body. For the \emph{Yellowfish} ASV, the force $F_u$ and torque $\tau_r$, expressed in the body-fixed frame \{b\}, can be found as:
\begin{equation}
    F_u = T(\delta_L) + T(\delta_R), \quad
    \tau_r = \frac{d}{2}\left(T(\delta_L) - T(\delta_R)\right),
    \label{eq_F_u_tau_r}
\end{equation}
where $d$ is the distance between propellers, and $T(\delta_i)$, with $i\in\{L,R\}$, represents the thrust dynamic model for any of the identical propellers.

Attending to the previous discussion, and to previous modelling in the literature, in this manuscript we have considered different models:
\begin{itemize}
    \item Function $T(\delta_i)$ is a static quadratic function of the PWM signal applied to the propeller. Considering the thrust quadratic with the RPM is common in the literature (see, for instance, \cite{simetti2022control}). Since the RPM is not a manipulable input, we seek a model that takes the PWM signal as input, assuming that this quadratic relation is preserved.
    \item Function $T(\delta_i)$ can be approximated by a first-order dynamic system. Again, this sort of model has been used in the literature (see, for example, \cite{blanke2000dynamic}). 
\end{itemize}

The next subsections present the aforementioned models. To simplify the notation, and assuming that both propellers are identical and thus share the same model, the subscripts $i\in\{L,R\}$ will be removed.

\subsubsection{Second-Order Static Model}\label{sec:static_model}

The thrust $T(\delta(k))$ is assumed to be described as:
\begin{equation}\label{eqThrust_2nd}
    T(\delta(k)) = \left\{\begin{array}{rl}
        a_{f}\delta(k)^2+b_{f}\delta(k),& \delta(k) \geq 0^+  \\
        a_{r}\delta(k)^2+b_{r}\delta(k),& \delta(k) < 0^-,
    \end{array}\right.
\end{equation}
where $\delta(k)\in[-1,1]$ is the normalised PWM signal, and $a_{j}$, $b_{j}\in\mathbb{R}$ with $j\in\{f,r\}$ are unknown scalars, where subindices $f$ and $r$ stand for forward and reverse motion, respectively. Equation \eqref{eqThrust_2nd} establishes four regions of operations attending to the rotation of the left and right propellers. The $(f,f)$ region will occur when the ASV moves forward in a straight line or with low-curvature trajectories. High-curvature trajectories are defined by the $(f,r)$ and $(r,f)$ regions. Finally, the reverse motion will be attained with the $(r,r)$ region.

As the primary objective of this identification process refers to the motion control of the \emph{Yellowfish} ASV, the latter configuration will not be employed during the identification experiments, as it does not align with the intended use-case scenario.

\begin{remark}
    The T200 Thrusters present a non-symmetric dead zone around $\delta=0$, which has been experimentally identified in \cite{morel2022modelling}. To account for it, \eqref{eqThrust_2nd} could be easily modified as:
\end{remark}

\begin{equation*}
    T(\delta(k)) = \left\{\begin{array}{rl}
        a_{f}\left(\delta(k)-\delta_{f}\right)^2+b_{f}\left(\delta(k)-\delta_{f}\right),& \delta(k) \geq \delta_f  \\
        a_{r}\left(\delta(k)-\delta_{r}\right)^2+b_{r}\left(\delta(k)-\delta_{r}\right),& \delta(k) \leq \delta_r,\\
        0,& \textrm{otherwise}
    \end{array}\right.
\end{equation*}
where $\delta_{f}>0$ and $\delta_{r}<0$ are the identified dead zone limits. The identification methods presented in the rest of this paper can be applied considering the previous model for the thrusters, or the one in~\eqref{eqThrust_2nd}. However, as the identified constants are very small, and trying to keep the notation as simple as possible to keep the focus on the identification, the rest of the document does not consider the dead zone.

\subsubsection{First-Order Dynamic Model}
Both thrusts are assumed to be described by a first-order model, this is:
\begin{equation}\label{eq_prop_dyn_model}    
T(k)= \alpha T(k-1)+\beta T_{st}(k-1)
\end{equation}
where $\alpha,\beta\in\mathbb{R}$ are constants to identify, and:
\begin{equation}\label{eqThrust_Poly_v2}
    T_{st}(\delta(k)) = \left\{\begin{array}{rl}
        a_{f}\delta(k)^2+b_{f}\delta(k),& \delta(k) \geq 0^+ \\
        a_{r}\delta(k)^2+b_{r}\delta(k),& \delta(k) \leq 0^-,
    \end{array}\right.
\end{equation}
represents the static input thrust as a quadratic function of the PWM normalised signal.

%%%%%%%%%%%%%%%%%%%%%%%%%%%%%%%%%%%%%%%%%%%%%%%%%
\section{Identification Procedure}
\label{sec:sample:Identification}
%%%%%%%%%%%%%%%%%%%%%%%%%%%%%%%%%%%%%%%%%%%%%%%%%
For the vehicle identification process, a specific procedure is followed as described in this section. First, the available data are prepared to meet the necessary requirements established in the previous section. Subsequently, the input gain for static propeller modelling is obtained. Finally, the equations are derived to identify the input gain, taking into account a dynamic model. This section provides a detailed explanation of the identification process for both propeller modelling.

According to Section \ref{sec:sample:ProblemStatement}, the problem to be addressed in the identification procedure is to find the input gains $G_u$, $G_v$, and $G_r$, whose mathematical expressions are defined in \eqref{eq_IG_definition2}, using the available experimental information.

\subsection{Data Preparation}\label{sec:data_prep}

The procedure to transform raw experimental data into the values necessary for identification is detailed. 
Initially, as data are collected at varying sampling times, the longest duration, in this case the GNSS data with a sampling time of $200$ ms, is used as a reference for the remainder. The input data are adjusted to use only previous points, with a cubic spline extrapolation applied to the latitude and longitude data from the GNSS. For the heading, which has a 20ms sampling time, an 8th-degree spline extrapolation is employed. For PWM data, a third-order spline extrapolation is also used. Please note that this synchronisation with spline extrapolations has been done carefully to keep causality, which can be lost with interpolation methods.

The values required for identification include the linear velocities $\{u,v\}$ and the yaw rate $r$, as well as the values of $\delta_L$ and $\delta_R$. However, instead of using individual PWM signals $\delta_{L}$ and $\delta_{R}$ as inputs, this manuscript proposes to use the PWM mean or $\bar{\delta}\equiv\frac{\delta_{L}+\delta_{R}}{2}$ and the PWM difference or $\Delta\delta\equiv \delta_{L}-\delta_{R}$. The reason is of a practical nature. In order to gather the necessary experimental data for the identification of the ASV, a series of experiments were conducted using a radio controller to issue commands to the propellers. Within this remote control system, two joysticks are configured to manipulate both the PWM mean and the PWM difference. It is not possible, with the radio controller, to choose the PWM values of each thruster independently. Furthermore, the PWM mean and difference signals are adequate to excite different kinds of dynamics, as will be explained later.

The process of obtaining the linear velocities $\{u,v\}$ is based on the measured positions and heading. This transformation also involves changing the reference frame from inertial coordinates to body-referenced coordinates. Detailed equations for this transformation are provided in~\cite{morel2022modelling}. Finally, these values are also filtered using a Savitzky-Golay filter. The yaw rate \(r\) is obtained by the derivative of $\psi$ and then filtered using a Savitzky-Golay filter~\cite{savitzky_golay}. Table~\ref{tab:identification_data} details the measured data from which these necessary values are derived.

\begin{table}[]
\centering
\caption{Required Data for Vessel Identification and Their Sources}
\label{tab:identification_data}
\resizebox{\columnwidth}{!}{%
\begin{tabular}{lll}
\toprule
\textbf{Data for Identification} & \textbf{Required Data} & \textbf{Description}                                                                                                                        \\ \hline
$u, v$                             & $x,y,\psi$                       & \begin{tabular}[c]{@{}l@{}}These derivatives obtain body-referenced velocities, \\ requiring heading.\end{tabular}                 \\
$\bar{\delta}$, $\Delta\delta $              & $\delta_L$, $\delta_R$            & \begin{tabular}[c]{@{}l@{}}The left and right thruster values are transformed \\ to obtain their mean and difference.\end{tabular} \\
$r$                                 &  $\psi$                             & \begin{tabular}[c]{@{}l@{}}The yaw rate is obtained by the derivative of $\psi$.\end{tabular}                     \\ \toprule
\end{tabular}%
}
\end{table}

\subsection{Input-Gain Identification for Static Propeller Modelling}

The left and right thrusts, $T(\delta(k))$, are assumed to be described as in~\eqref{eqThrust_2nd}. According to \eqref{eq_F_u_tau_r}, now using the PWM mean $\bar{\delta}$ and  difference $\Delta \delta$ values, the force $F_u$ and the torque $\tau_r$ can be written as follows:

\begin{align}
    F_u(\bar{\delta}(k),\Delta\delta(k)) &= (a_{j}^{L}+a_{j}^{R})\left(\bar{\delta}(k)^2+\frac{\Delta\delta(k)^2}{4}\right)+(a_{j}^{L}-a_{j}^{R})\bar{\delta}(k)\Delta\delta(k)\nonumber\\
     &\quad + (b_{j}^{L}+b_{j}^{R})\bar{\delta}(k)+(b_{j}^{L}-b_{j}^{R})\frac{\Delta\delta(k)}{2},\label{eq_F_u_2}\\
    \tau_r(\bar{\delta}(k),\Delta\delta(k)) &= \frac{d}{2}\Bigg((a_{j}^{L}-a_{j}^{R})\left(\bar{\delta}(k)^2+\frac{\Delta\delta(k)^2}{4}\right)+(a_{j}^{L}+a_{j}^{R})\bar{\delta}(k)\Delta\delta(k)\nonumber\\
     &\quad + (b_{j}^{L}-b_{j}^{R})\bar{\delta}(k)+(b_{j}^{L}+b_{j}^{R})\frac{\Delta\delta(k)}{2}\Bigg),\label{eq_tau_r_2}
\end{align}
where $a_{j}^{i},b_{j}^{i}$ with $i\in\{L,R\},j\in\{f,r\}$ are the forward and reverse models to be identified for the propellers. 

Now, for a better excitation of the surge dynamics, only those instants in the dataset in which the propellers work in the region $(f,f)$ will be considered for identification. In addition, the external disturbances are considered constant for identification, this is, $c_u(k)=\bar{c}_{u},\;\forall k$. 
\begin{prop}
    Provided that $\bar{\delta}(k)$ and $\Delta\delta(k)$ are known for every $k$, and the velocities $u(k),v(k),r(k)$ can be obtained from the position data as explained in Section~\ref{sec:data_prep}, the evolution of the surge at instant $k$ for a quadratic static model for the propeller and operating in the region $(f,f)$ satisfies:
\end{prop}

\begin{equation}\label{eq:linear_map_surge}
    A_{u}(k)X_{u}=b_{u}(k),
\end{equation}
where $X_u$ represents a vector of unknown coefficients to identify:
\begin{equation}
    X_u = h[P_{u}^{(1,1)}~~2Q_{u}^{(2,3)}~~Q_{u}^{(3,3)}~~R_{u}^{(1)}~~\bar{c}_u~~2\bm{M}^{-1}(1,1)a_f~~2\bm{M}^{-1}(1,1)b_f]^{T},
\end{equation}
and $A_u(k),b_u(k)$ are a vector and a scalar, respectively, that can be computed from the dataset:
\small
\begin{align}
    A_u(k) &= \left[u(k)|u(k)|~~v(k)r(k)~~r(k)^2~~u(k)~~1~~\left(\bar{\delta}(k)^2+\frac{\Delta{\delta}(k)^2}{4}\right)~~\bar{\delta}(k)\right],   \\
    b_u(k) &= u(k+1) - u(k).
\end{align}
\normalsize

\begin{proof}
     While in the region $(f,f)$, the force in~\eqref{eq_F_u_2} can be simplified as:
\begin{equation*}
    F_u(\bar{\delta}(k),\Delta\delta(k)) =  2a_{f}\left(\bar{\delta}(k)^2+\frac{\Delta\delta(k)^2}{4}\right)+ 2b_{f}\bar{\delta}(k).
\end{equation*}

Subsequently, the dynamics of the surge for a quadratic static model for the propeller results in:
\begin{align*}
    u(k+1) &=u(k)+2h\bm{M}^{-1}(1,1)a_{f}\left(\bar{\delta}(k)^2 + \frac{\Delta{\delta}(k)^2}{4}\right) \nonumber \\ 
    &\quad+2h\bm{M}^{-1}(1,1)b_{f}\bar{\delta}(k) + hP_{u}^{(1,1)}u(k)|u(k)|  \nonumber\\ 
&\quad+2hQ_{u}^{(2,3)}v(k)r(k)+hQ_{u}^{(3,3)}r(k)^2+hR_{u}^{(1)}u(k)+hc_{u}(k).
\end{align*}

Since the mean and difference PWM normalised signals are known, and the velocities can be obtained from the dataset, the previous equation reveals that the unknown coefficients appear in a linear fashion, as in \eqref{eq:linear_map_surge}. 
\end{proof}

\begin{remark}
    The experiments were carried out on an artificial lake that had minimal influence from currents and waves. Concerning the wind, the experiments took place on relatively calm days. However, it is not possible to completely neglect the effect of the wind on lumped disturbances. To partially compensate for that, the identification procedure tries to capture a constant disturbance term $\bar{c}_{u}$ with the intention of modelling a bias or persistent force in some direction.
\end{remark}

Vector $X_u$ can be approximated, for instance, via least-squares. As a result, the surge input gain $G_u(k)$ in~\eqref{eq_IG_definition2} can be computed as a function of the mean and difference PWM normalised values:
\begin{equation}
    G_u(k) = X_u(6)\left(\bar{\delta}(k)^2+\frac{\Delta{\delta}(k)^2}{4}\right)+X_u(7)\bar{\delta}(k).
    \label{eq_IG:u}
\end{equation}

Concerning the yaw and sway dynamics, now the data belonging to regions $(f,f),(f,r),(r,f)$ will be used for identification. Despite the fact that sway and, particularly, yaw dynamics are better excited in regions $(f,r),(r,f)$ than in region $(f,f)$, the \emph{Yellowfish} ASV is mainly designed to conduct path-following missions and, therefore, high-curvature trajectories will not be predominant.
\begin{prop} \label{prop:sway_yaw}
    Provided that $\bar{\delta}(k)$ and $\Delta\delta(k)$ are known for every $k$, and the velocities $u(k),v(k),r(k)$ can be obtained from the position data, the evolution of the sway and yaw at instant $k$ for a quadratic static model for the propeller satisfy:
\begin{align}\label{eq:linear_map_sway_yaw}
    A_{v}(k)X_{v}&=b_{v}(k),\\
    A_{r}(k)X_{r}&=b_{r}(k),
\end{align}
where $X_v$ and $X_r$ are the vectors of unknown coefficients to identify:
\begin{align}
    X_v &= h\Big[P_{v}^{(2,2)}~~P_{v}^{(2,3)}~~P_{v}^{(3,2)}~~P_{v}^{(3,3)} ~~2Q_v^{(1,2)}~~2Q_v^{(1,3)}~~R_v^{(2)}~~R_v^{(3)}~~\Bar{c}_v \nonumber \\ 
    &\quad~~\bm{M}^{-1}(2,3)\frac{d}{2}(a_{j}^{L}-a_{j}^{R})~~\bm{M}^{-1}(2,3)\frac{d}{2}(a_{j}^{L}+a_{j}^{R})\nonumber \\
    &\quad~~\bm{M}^{-1}(2,3)\frac{d}{2}(b_{j}^{L}-b_{j}^{R})~~\bm{M}^{-1}(2,3)\frac{d}{2}(b_{j}^{L}+b_{j}^{R})\Big]^T,\\
    X_r &= h\big[P_{r}^{(2,2)}~~P_{r}^{(2,3)}~~P_{r}^{(3,2)}~~P_{r}^{(3,3)} ~~2Q_r^{(1,2)}~~2Q_r^{(1,3)}~~R_r^{(2)}~~R_r^{(3)}~~\Bar{c}_r \nonumber \\ 
    &\quad~~\bm{M}^{-1}(3,3)\frac{d}{2}(a_{j}^{L}-a_{j}^{R})~~\bm{M}^{-1}(3,3)\frac{d}{2}(a_{j}^{L}+a_{j}^{R})\nonumber \\
    &\quad~~\bm{M}^{-1}(3,3)\frac{d}{2}(b_{j}^{L}-b_{j}^{R})~~\bm{M}^{-1}(3,3)\frac{d}{2}(b_{j}^{L}+b_{j}^{R})\big]^T,
\end{align}
and $A_v(k),A_r(k),b_v(k),b_r(k)$ are two pairs of vectors and scalars, respectively, that can be computed from the dataset. 
\end{prop}

\begin{proof}
    To follow the step-by-step process of Proposition~\ref{prop:sway_yaw}, the reader is referred to~\ref{apendix_static_swayyaw}.
\end{proof}

As before, vectors $X_v,X_r$ can be approximated, for instance, via least-squares. Then, sway and yaw input gains can be computed from the mean and difference PWM normalised values:

\begin{equation}
   G_p(k) =\left\{\begin{array}{ll}
   X_p(11)\bar{\delta}(k)\Delta\delta(k) + X_p(13)\frac{\Delta\delta(k)}{2}, & (f,f) \\
   \begin{array}{r}
   X_p(10)\left(\bar{\delta}(k)^2+\frac{\Delta\delta(k)^2}{4}\right) +X_p(11)\bar{\delta}(k)\Delta\delta(k)  \\
   +X_p(12)\bar{\delta}(k) +X_p(13)\frac{\Delta\delta(k)}{2}\end{array}, & (f,r)\\
   \begin{array}{r}
   -X_p(10)\left(\bar{\delta}(k)^2+\frac{\Delta\delta(k)^2}{4}\right) +X_p(11)\bar{\delta}(k)\Delta\delta(k)  \\
   -X_p(12)\bar{\delta}(k) +X_p(13)\frac{\Delta\delta(k)}{2}\end{array}, & (r,f)
   \end{array}\right.
   \label{eq:IG_vr}
\end{equation}
where $p\in\{v,r\}$.

\begin{remark}
    Please note that, even though it is not the objective of the method, the previous procedure allows for a reconstruction of the lumped disturbances as a function of the velocities. This could be interesting when designing a feedforward controller that could partially compensate for the effect of such disturbances.
\end{remark}

\subsection{Input-Gain Identification for Dynamic Propeller Modelling}
\label{sec:dynamic_identification}

In the process of identifying the input gain considering the dynamic propeller modelling, \eqref{eq_prop_dyn_model} is used for the first-order dynamic model. Following this, the procedure to obtain the regressors bears a resemblance to the approach used to identify the input gain for static propeller modelling and, hence, some developments will be simplified. The following aims to highlight the specific adaptations necessary for dynamic analysis.

\begin{prop} \label{prop:dynamic}
   Provided that $\bar{\delta}(k)$ and $\Delta\delta(k)$ are known for every $k$, and the velocities $u(k),v(k),r(k)$ can be obtained from the position data, the evolution of the velocities at instant $k$ for a quadratic dynamic model for the propeller satisfies: 
\begin{align}
    A_{dyn,u}(k)X_{dyn,u}&=b_{dyn,u}(k), \label{eq:surgedyn}\\ 
    A_{dyn,v}(k)X_{dyn,v}&=b_{dyn,v}(k), \label{eq:swaydyn}\\
    A_{dyn,r}(k)X_{dyn,r}&=b_{dyn,r}(k), \label{eq:yawdyn}
\end{align}
where $X_{dyn,u}$, $X_{dyn,v}$ and $X_{dyn,r}$ are the vectors of unknown coefficients to identify:
\begin{align}
    X_{dyn,u} &= \Bigg[\left( \alpha+hR_{u}^{(1)} \right)~~-\alpha hP_{u}^{(1,1)}~~-\alpha hQ_{u}^{(2,3)}~~-\alpha hQ_{u}^{(3,3)}\nonumber \\ &\quad~~- \alpha\left( 1+hR_{u}^{(1)} \right)~~hP_{u}^{(1,1)}~~hQ_{u}^{(2,3)}~~hQ_{u}^{(3,3)}\nonumber \\
        &\quad~~(h-\alpha)\bar{c_{u}}~~2h\beta\bm{M}^{-1}(1,1)a_f~~2h\beta\bm{M}^{-1}(1,1) b_f\Bigg]^T,
\end{align}
\begin{align}
    X_{dyn,v} &= \Bigg[\left( \alpha+hR_{v}^{(2)} \right)~~-\alpha hP_{v}^{(2,2)}~~-\alpha hP_{v}^{(2,3)}~~-\alpha hP_{v}^{(3,2)}~~-\alpha hP_{v}^{(3,3)}\nonumber \\ &\quad~~-2\alpha hQ_{v}^{(1,2)}~~-2\alpha hQ_{v}^{(1,3)}~~-\alpha \left(1+hR_{v}^{(2)}\right)~~-\alpha hR_{v}^{(3)}~~ hP_{v}^{(2,2)}\nonumber \\ 
    &\quad~~hP_{v}^{(2,3)}~~hP_{v}^{(3,2)}~~hP_{v}^{(3,3)}~~2hQ_{v}^{(1,2)}~~2hQ_{v}^{(2,2)}~~hR_{v}^{(2)}~~hR_{v}^{(3)}~~(h-\frac{d}{2}\alpha)\bar{c_{v}}\nonumber \\ &\quad~~\frac{d}{2}h\beta\bm{M}^{-1}(2,3)(a_j^L-a_j^R)~~\frac{d}{2}h\beta\bm{M}^{-1}(2,3)(a_j^L+a_j^R)\nonumber \\ &\quad~~\frac{d}{2}h\beta\bm{M}^{-1}(2,3)(b_j^L-b_j^R)~~\frac{d}{2}h\beta\bm{M}^{-1}(2,3)(b_j^L+b_j^R)\Bigg]^T, \\
     X_{dyn,r} &= \Bigg[\left( \alpha+hR_{r}^{(3)} \right)~~-\alpha hP_{r}^{(2,2)}~~-\alpha hP_{r}^{(2,3)}~~-\alpha hP_{r}^{(3,2)}~~-\alpha hP_{r}^{(3,3)}\nonumber \\ &\quad~~-2\alpha hQ_{r}^{(1,2)}~~-2\alpha hQ_{r}^{(1,3)}~~-\alpha hR_{r}^{(2)}~~-\alpha\left(1+hR_{r}^{(3)}\right)~~hP_{r}^{(2,2)}\nonumber \\ &\quad~~hP_{r}^{(2,3)}~~hP_{r}^{(3,2)}~~hP_{r}^{(3,3)}~~2hQ_{r}^{(1,2)}~~2hQ_{r}^{(2,2)}~~hR_{r}^{(2)}~~hR_{r}^{(3)}~~(h-\frac{d}{2}\alpha)\bar{c_{r}}\nonumber \\ &\quad~~\frac{d}{2}h\beta\bm{M}^{-1}(3,3)(a_j^L-a_j^R)~~\frac{d}{2}h\beta\bm{M}^{-1}(3,3)(a_j^L+a_j^R)\nonumber \\ &\quad~~\frac{d}{2}h\beta\bm{M}^{-1}(3,3)(b_j^L-b_j^R)~~\frac{d}{2}h\beta\bm{M}^{-1}(3,3)(b_j^L+b_j^R)\Bigg]^T.
\end{align}
\end{prop}

\begin{proof}
    For a detailed breakdown of Proposition~\ref{prop:dynamic}, consult~\ref{apendix_dynamic}.
\end{proof}

In the identification of the input gain considering the static model of the propellers, the regressors for each component of the velocity contained a different set of parameters and therefore the identification could be tackled in three independent steps. This is not exactly the case here. Observing, for instance, components 1 and 5 of the regressor $X_{dyn,u}$, it is easy to see that there are 2 unknowns $\left(\alpha,R_u^{(1)}\right)$, since $h$ is known. Hence, the value of $\alpha$ can be identified using just the regressor of the surge. However, from the components 1 and 8 of the regressor $X_{dyn,v}$, one can identify $\left(\alpha,R_v^{(3)}\right)$, obtaining, possibly, a different value of $\alpha$. Something similar happens with the regressor of the yaw. This problem does not occur with $\beta$, as it appears to always be coupled with other parameters that are different for each component of the velocity.

To solve this issue, there are different options. One possible option would consist in considering different dynamic models for each component, in other words, a different $\alpha_j$ for each $j\in\{u,v,r\}$. Although possible from a numerical perspective, this option does not look fully consistent with the underlying physics. Another option, the one implemented here, is to identify a unique $\alpha$ valid for all velocities. This implies solving an overdetermined system of equations, presented in~\eqref{eq_alpha}, since there are more equations than unknowns in the regressors of the three velocities.

\begin{equation}
\begin{cases}
\label{eq_alpha}
X_{\text{dyn},u}(1) &= \alpha + hR_u^{(1)} \\
X_{\text{dyn},u}(5) &= -\alpha(1 +h R_u^{(1)}) \\
X_{\text{dyn},v}(1) &= \alpha + hR_v^{(2)} \\
X_{\text{dyn},v}(8) &= -\alpha(1 + hR_v^{(2)}) \\
X_{\text{dyn},r}(1) &= \alpha + hR_r^{(3)} \\
X_{\text{dyn},r}(9) &= -\alpha(1 + hR_r^{(3)})
\end{cases}
\end{equation}

\begin{remark}
    Please note that, with the proposed identification procedure considering dynamic propeller modelling, the value of some parameters related to the lumped disturbances, such as $P_u^{(1,1)}$, could have two different values. Although this is, of course, impossible, it is not particularly important since we are interested in identifying the input gain, which will be uniquely determined.
\end{remark}

Once the parameters are identified, the input gain can be obtained. From the definition of the input gain in~\eqref{eq_IG_definition2}, using the force and torque in~\eqref{eq_F_u_tau_r}, and the first-order dynamic model of the thrusts in~\eqref{eq_prop_dyn_model}, the following dynamic model for the input gain is formulated as:
\begin{align} \label{eq:dynk}
    \bm{G_{dyn}}(k) &= h\bm{M}^{-1}\alpha\bm{\tau}(k-1) + h\bm{M}^{-1}\beta\bm{\tau_{st}}(k-1),
\end{align}
where $\bm{G_{dyn}}(k)=[G_{dyn,u}(k), G_{dyn,v}(k), G_{dyn,r}(k)]^T$ and $\bm{\tau_{st}}(k)$ is made up of the force $F_u$ and torque $\tau_r$ given by~\eqref{eq_F_u_2} and~\eqref{eq_tau_r_2}, which turn out to be the static model equations for the force and torque. 
Subsequently, we can group the expression $h\bm{M}^{-1}\bm{\tau}(k-1)$ as $\bm{G_{dyn}}(k-1)$. Taking advantage of this simplification, \eqref{eq:dynk} becomes:
 \begin{align}
     \bm{G_{dyn}}(k) &= \alpha\bm{G_{dyn}}(k-1) + h\bm{M}^{-1}\beta\bm{\tau_{st}}(k-1),
 \end{align}
which shows that the input gain allows now for a first-order dynamic model.

This sets the groundwork for applying the equation across all velocities, culminating in a comprehensive application of the dynamic model for input gain estimation as shown:
\begin{align} \label{eq_IG_dyn_u}
     G_{dyn,u}(k) &= \alpha G_{dyn,u}(k-1) + X_{dyn,u}(10)\left(\bar{\delta}(k-1)^2+\frac{\Delta{\delta}(k-1)^2}{4}\right)\nonumber \\ &\quad+X_{dyn,u}(11)\bar{\delta}(k-1),
\end{align}
\small
\begin{equation} \label{eq_IG_dyn_vr}
   G_{dyn,p}(k) =\left\{\begin{array}{lr}
   \alpha G_{dyn,p}(k-1) + X_{dyn,p}(19)\bar{\delta}(k-1)\Delta{\delta(k-1)} 
   \\+ X_{dyn,p}(21)\frac{\Delta{\delta(k-1)}}{2}, & (f,f) \\
   \begin{array}{l}
   \alpha G_{dyn,p}(k-1) \pm X_{dyn,p}(18)\left(\bar{\delta}(k-1)^2+\frac{\Delta{\delta}(k-1)^2}{4}\right) \\ + X_{dyn,p}(19)\bar{\delta}(k-1)\Delta{\delta(k-1)} \pm X_{dyn,p}(20)\bar{\delta}(k-1) \\ + X_{dyn,p}(21)\frac{\Delta{\delta(k-1)}}{2}\end{array}, & (f,r)\\
   \begin{array}{l}
   \alpha G_{dyn,p}(k-1) \pm X_{dyn,p}(18)\left(\bar{\delta}(k-1)^2+\frac{\Delta{\delta}(k-1)^2}{4}\right) \\ + X_{dyn,p}(19)\bar{\delta}(k-1)\Delta{\delta(k-1)} \pm X_{dyn,p}(20)\bar{\delta}(k-1) \\ + X_{dyn,p}(21)\frac{\Delta{\delta(k-1)}}{2}\end{array}, & (r,f)
   \end{array}\right.
\end{equation}
\normalsize
where $p\in\{v,r\}$.

Being a first-order model, the input gain cannot be computed only with the information of the PWM signals but also requires an initial condition. This initial condition is generally unknown. However, we can take advantage of the following property.

\begin{proper}
    Consider the evolution of the input gain given in~\eqref{eq_IG_dyn_u}-\eqref{eq_IG_dyn_vr}. Then, provided that $\alpha<1$, the value of the input gain at any instant $k'$ sufficiently large mainly depends on the value of the PWM signals. In other words, the influence of the initial condition vanishes with time.
\end{proper}

\begin{proof}
    The proof is straightforward since the input gain is a first-order system with a transition matrix $\Phi(k)=\alpha^k$. Provided that $\alpha<1$, the influence of the initial conditions vanishes as $k$ increases.
\end{proof}

\begin{remark} \label{remark:alpha}
    The fact that $\alpha<1$ could be inferred from the definition of the dynamic model for the thrust in~\eqref{eq_prop_dyn_model} since otherwise the thrust could grow to infinity with zero input. However, without imposition, a suitable value of $\alpha$ is obtained with the proposed identification procedure, as shown in \eqref{sec:sample:Results}.
\end{remark}

%%%%%%%%%%%%%%%%%%%%%%%%%%%%%%%%%%%%%%%%%%%%%%%%%
\section{Identification Experiments, Results and Discussion}
\label{sec:sample:Results}
%%%%%%%%%%%%%%%%%%%%%%%%%%%%%%%%%%%%%%%%%%%%%%%%%
In this section, the conducted experiments and the obtained results are presented and discussed. First, a detailed explanation of the experimental setup and the identification metrics used is provided. The identification results are then displayed using the second-order static model. The section concludes with the results obtained from the first-order dynamic model, including a discussion comparing both analysed models.

\subsection{Experimental Setup and Identification Metrics}

The \emph{Yellowfish} ASV was programmed to carry out a wide range of trajectories, divided into two main categories: forward and circular movements. The aim was to collect significant experimental data to excitate the surge and yaw dynamics. In the case of forward movements, manoeuvres were performed using constant and variable forward velocities, as well as zigzag trajectories. Similarly, for circular movements, manoeuvres were designed to excite yaw dynamics at both constant and variable rates, also integrating movements with constant and variable forward velocity into the turns. 

This extensive set of manoeuvres was intended to thoroughly evaluate the ASV performance in executing simple linear paths and complex navigational patterns. These manoeuvres produced a detailed dataset crucial for analysing and validating the navigational model, accurately reflecting the ASV real-world performance under a variety of controlled conditions.

The dataset was then strategically divided into two subsets: one devoted to the identification of the surge dynamics, whose points correspond to those in which forward movements were performed, where both propellers are in the $\{f,f\}$ zone; and another one for the identification of the yaw and sway dynamics, including those points that correspond to circular movements, considering the $\{f,r\}$ and $\{r,f\}$ zones and always disregarding any point in the $\{r,r\}$ zone, as it will not be analysed in this identification process. 
Then, each subset is further divided into a training set and a validation set. The complete dataset, comprising 16,025 data points, with a higher concentration of data within the $\{f,r\}$ and $\{r,f\}$ zones compared to the $\{f,f\}$ zone, divided into 127 trajectory segments, resulting in 53.42 minutes of experiments, is available in~\cite{Morel2024}. 

\begin{figure}[h]
    \begin{centering}
    \includegraphics[angle=0,width=0.62\textwidth]{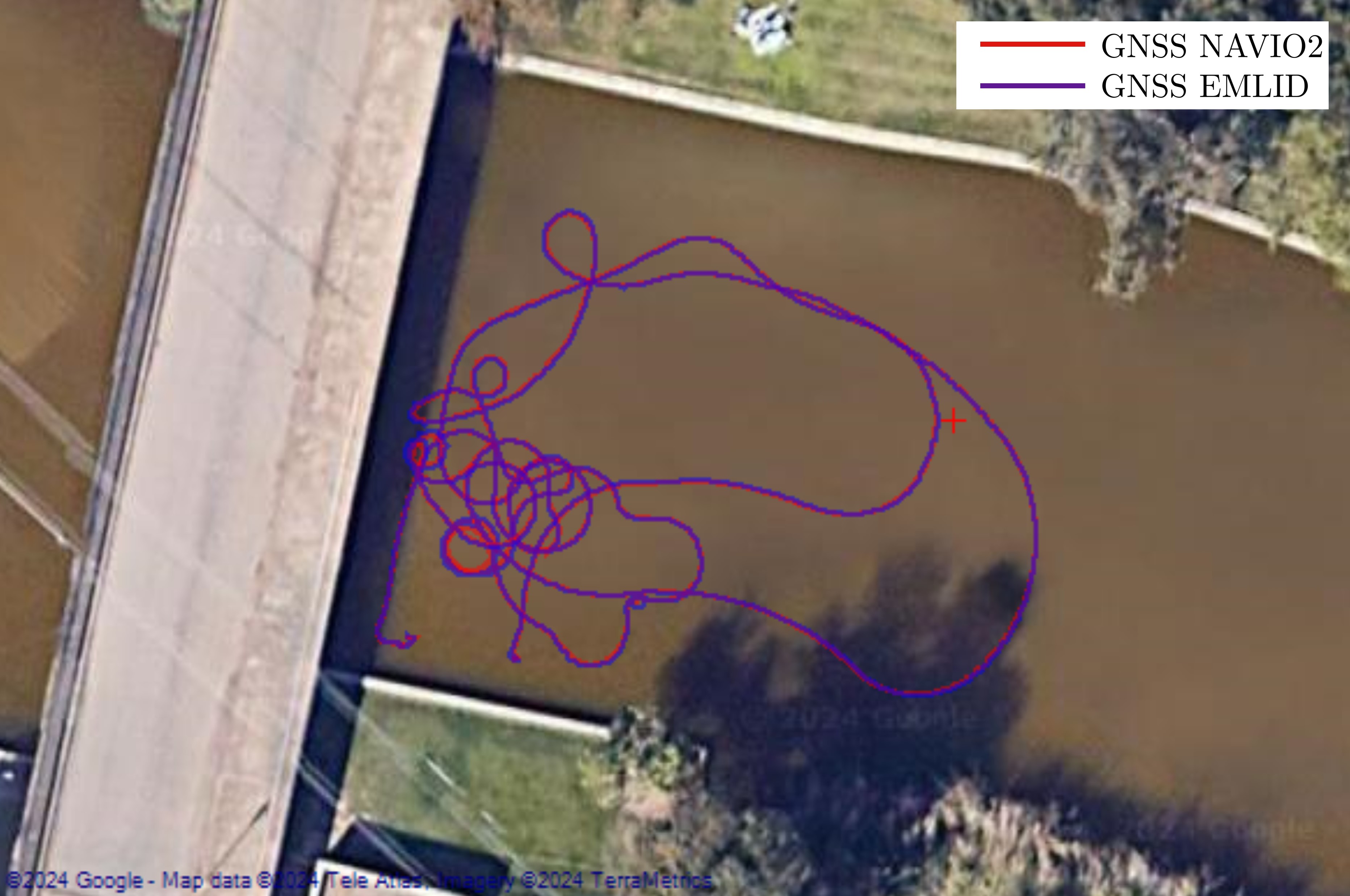}
    \par\end{centering}
    \caption{The Experiment Field and ASV Path}
    \label{fig_ASVpath}
\end{figure}

The way in which the dataset is split into training and validation subsets is done in two different ways. First, all data points are randomly divided, according to a given percentage, into these two sets. The second option randomly divides the trajectory segments, so that some of them belong to the training dataset, and others to the validation dataset. Since the segments have a fixed number of data points, this partition will not achieve an exact division of the dataset in terms of a given percentage but will try to approximate it as much as possible. Table~\ref{tab:separation_method} presents the separation methods used, along with a description of each; Figure~\ref{fig_ASVpath} shows the field in which the experiments were conducted. It depicts a calm lake setting. Two paths are presented: one generated by the internal Navio2 GNSS and the other by the external GNSS, which is the one used for its higher satellite precision.

\begin{table}[b]
\centering
\caption{Description of Data Separation Methods for Training Data}
\label{tab:separation_method}
\resizebox{\columnwidth}{!}{%
\begin{tabular}{ll}
\toprule
\multicolumn{1}{c}{\textbf{Separation Method}} & \multicolumn{1}{c}{\textbf{Description}}                                \\ \hline
By Segments                                    & Random grouping of the complete segments.                       \\
By Independent Points                          & Random and independent selection of the individual data points. \\ \toprule
\end{tabular}%
}
\end{table}

Regarding the values used for the validation analysis, the velocity values obtained by the data preparation process of the experimental results~($\bm{\nu}$) are the ground truth and are compared to the velocities predicted with the identified parameters~($\bm{\hat{\nu}})$.

Two well-known metrics are used for validation, the $R^2$, or coefficient of determination, and the mean absolute error (MAE): 

\begin{itemize}
\item The $R^2$ is a statistical measure that represents the proportion of variance in a dependent variable that is predictable from the independent variables. It provides an indication of how well the data points fit a model; the closer $R^2$ is to 1, the better the model explains the variability of the output in relation to the input variables. It is a measure commonly used in least-squares regressions~\cite{cameron1997r}. 
 
For the problem at hand, this metric is computed as follows:

\begin{equation}
    R_{j}^2 = 1 - \frac{\sum_{k=1}^n (\bm{\nu}_j(k) - \hat{\bm{\nu}}_j(k))^2}{\sum_{k=1}^n (\bm{\nu}_j(k) - \overline{\bm{\nu}}_j)^2},
    \label{eq:Rsquared}
\end{equation}
where $j\in\{u,v,r\}$ and $n$ is the number of points of the validation set.

\item The MAE is particularly effective for regressors that involve multiple parameters because it provides a straightforward measure of the average magnitude of errors without considering
their direction, which makes it highly suitable for evaluating the performance of models in predicting a wide range of outcomes~\cite{willmott2005advantages}. In this case, it is computed as:

\begin{equation}
    \text{MAE}_j = \frac{1}{n} \sum_{k=1}^n |\bm{\nu}_j(k) - \bm{\hat{\nu}}_j(k)|.
    \label{eq:MAE}
\end{equation}
\end{itemize}

\subsection{Identification Results for Static Propeller Modelling}

In the first analysis, the sensitivity of the identification results with respect to the partition of the data is studied. To do so, a fixed 70\% of the dataset is used for training and the 30\% left is used for validation. To analyse the sensitivity, 20 random partitions are made and, for each metric, the mean and Standard Deviation (SD) are computed. Table~\ref{tab:MAE_R2_varianza} presents the numerical results for the two partition methods.

From these data, two important conclusions can be drawn. First, the random partition in independent points renders slightly better results. But, most importantly, is the fact that the deviation is very small (much smaller with independent points). This allows us to conclude that the way in which the data are partitioned does not have an important impact on the predictive capabilities of the identified model. This conclusion invites us to avoid the repetitions of the partitions in the rest of the analysis.

Attending to the quantitive values observed, the identified models show a really accurate prediction, in surge, sway, and yaw, with a determination coefficient always higher than 0.98.

\begin{table}[h]
\centering
\caption{Comparison of MAE and $R^2$ Metrics for Different Data Separation Methods Using Second-Order Static Model}
\label{tab:MAE_R2_varianza}
\resizebox{\columnwidth}{!}{%
\begin{tabular}{lllllll}
\toprule
                 & \textbf{MAE}$_u$ & \textbf{SD}$_u$ & \textbf{MAE}$_v$ & \textbf{SD}$_v$ & \textbf{MAE}$_r$ & \textbf{SD}$_r$ \\ \cline{2-7} 
\textbf{By Segments} &  0.014377 & $2.9\times 10^{-3}$      &   0.011913 & $6.2\times 10^{-4}$     & 0.015895 & $2.1\times 10^{-3}$       \\
\textbf{By Independent Points}   &  \bm{$0.011794$} & \bm{$1.5\times 10^{-4}$}      &   \bm{$0.011870$} & \bm{$1.0\times 10^{-4}$}       &  \bm{$0.015757$} &\bm{$2.0\times 10^{-4}$}       \\ \hline
                 & \bm{$R^2_u$}  & \textbf{SD}$_u$ & \bm{$R^2_v$} & \textbf{SD}$_v$ & \bm{$R^2_r$} & \textbf{SD}$_r$ \\ \cline{2-7} 
\textbf{By Segments} & 0.995708 & $2.3\times 10^{-3}$   & 0.987277 & $1.8\times 10^{-3}$       &  0.997468 & $1.2\times 10^{-3}$      \\
\textbf{By Independent Points}  & \bm{$0.997658$} & \bm{$9.6\times 10^{-5}$}       &   \bm{$0.987452$} & \bm{$3.6\times 10^{-4}$}     & \bm{$0.997659$} & \bm{$7.4\times 10^{-5}$}    \\ \toprule
\end{tabular}%
}
\end{table}

In the next analysis, the sensitivity to the number of points in the training dataset is analysed. To do so, different percentages of points in the whole dataset are considered, with a constant keeping 30\% for the validation dataset. The results are shown in Tables~\ref{tab:MAE_percentages_static} and~\ref{tab:R2_percentages_static}.

With respect to the MAE, there is no significant change as the training dataset is reduced. This allows us to conclude that the performance of the model is consistent despite the reduction of data points.

Additionally, since a larger percentage of training might capture more variability and complexity, these results could suggest that the model is not overfitting to the training data. This is indicated by the fact that, in some cases, the MAE values for the validation set are lower than those for the training set, implying that the model is not overly fine-tuned to the specific details of the training data.

\begin{table}[h!]
\centering
\caption{MAE Values for Different Training Percentages Using Second-Order Static Model}
\label{tab:MAE_percentages_static}
\resizebox{\columnwidth}{!}{%
\begin{tabular}{cllllcllll}
\toprule
\multicolumn{1}{l}{\textbf{Training \%}} & \textbf{MAE}$_u$      & \textbf{MAE}$_v$      & \textbf{MAE}$_r$      &  & \multicolumn{1}{l}{\textbf{Validation \%}} & \textbf{MAE}$_u$      & \textbf{MAE}$_v$      & \textbf{MAE}$_r$      &  \\ \cline{1-4} \cline{6-9}
70                              & \bm{$0.012234$} & \bm{$0.012666$} & \bm{$0.013471$} &  & 30                                & \bm{$0.011092$} & 0.010599 & 0.020743 &  \\
60                              & 0.013645 & 0.012879 & 0.014153 &  & 30                                & 0.011370 & 0.010596 & 0.020164 &  \\
50                              & 0.016374 & 0.014090 & 0.014068 &  & 30                                & 0.011580 & \bm{$0.010361$} & \bm{$0.019268$}&  \\ 
\bottomrule
\end{tabular}
}
\end{table}

Comparing the results using $R^2$, a concordance with the obtained MAE values can be observed, allowing us to conclude that the model achieves a good level of prediction.

\begin{table}[h!]
\centering
\caption{$R^2$  Values for Different Training Percentages Using Second-Order Static Model}
\label{tab:R2_percentages_static}
\resizebox{\columnwidth}{!}{%
\begin{tabular}{cllllcllll}
\toprule
\textbf{Training \%} & \bm{$R^2_u$ }     & \bm{$R^2_v$}      & \bm{$R^2_r$ }     &  & \textbf{Validation \%} & \bm{$R^2_u$ }     & \bm{$R^2_v$}      &\bm{ $R^2_r$ }    &  \\ \cline{1-4} \cline{6-9}
70       & \bm{$0.998177$} & 0.984567 & \bm{$0.998602$} &  & 30         & \bm{$0.987774$} & 0.991418 & 0.992662 &  \\
60       & 0.996305 & \bm{$0.985519$} & 0.997190 &  & 30         & 0.987496 & 0.991008 & 0.992628 &  \\
50       & 0.988801 & 0.982918 & 0.997262 &  & 30         & 0.987273 & \bm{$0.991643$} & \bm{$0.993209$} &  \\ \toprule
\end{tabular}%
}
\end{table}

Finally, some graphical analysis are shown. Although the metrics provide very useful information, a graphical analysis allows us to understand the behaviour of the model, when it renders nice predictions, and when it fails. For this part, the dataset is partitioned in segments, this allowing us to observe dynamical behaviours.

As mentioned previously, the dataset was strategically divided into two subsets: one dedicated to surge dynamics and the other to yaw and sway dynamics. Consequently, the graphical representation of the surge results uses a different time window compared to those of sway and yaw, as they exhibit different movements.

\begin{figure}[h!]
    \centering
    \subfigure[]{\includegraphics[width=0.32\textwidth]{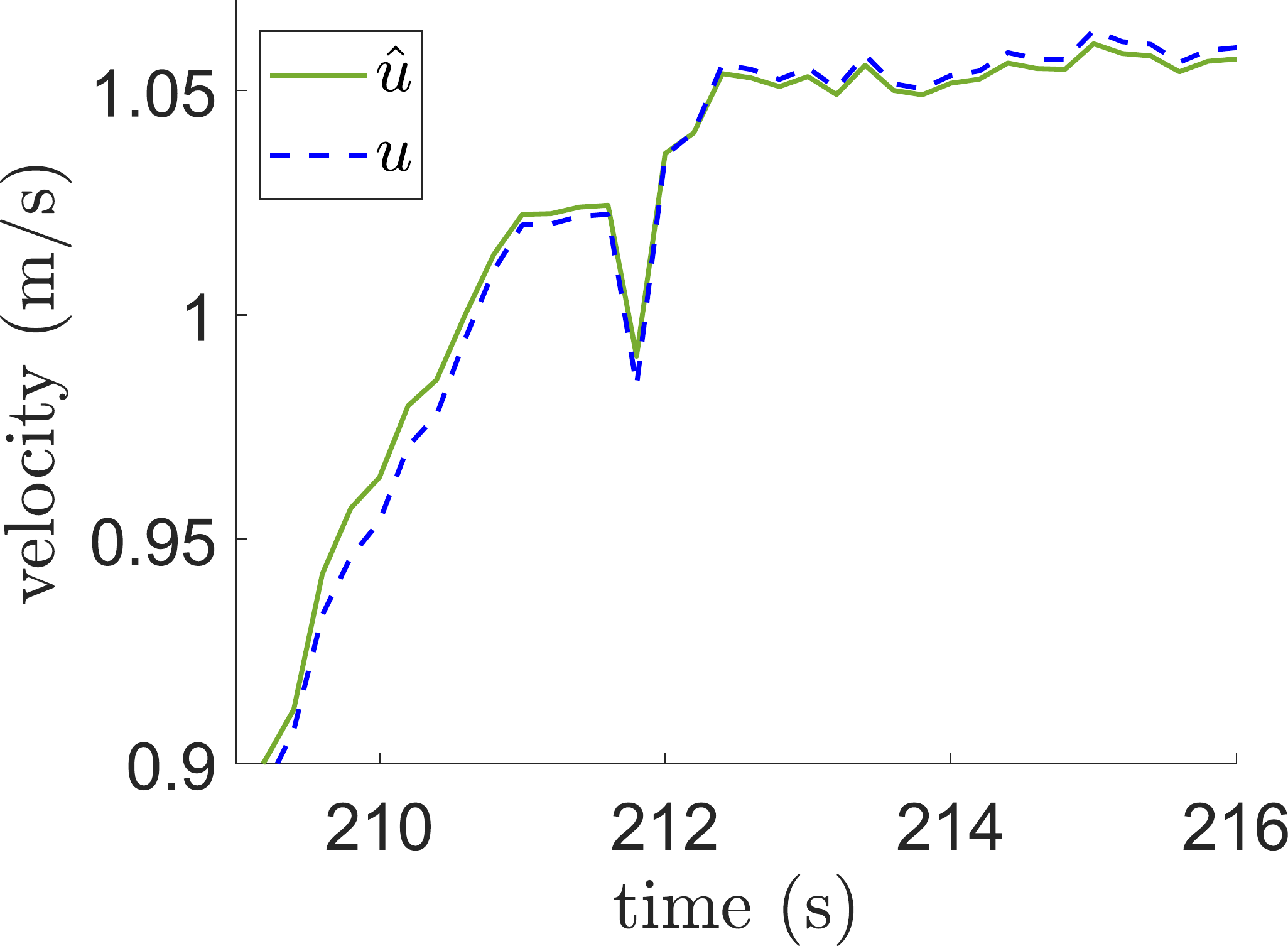}} 
    \subfigure[]{\includegraphics[width=0.32\textwidth]{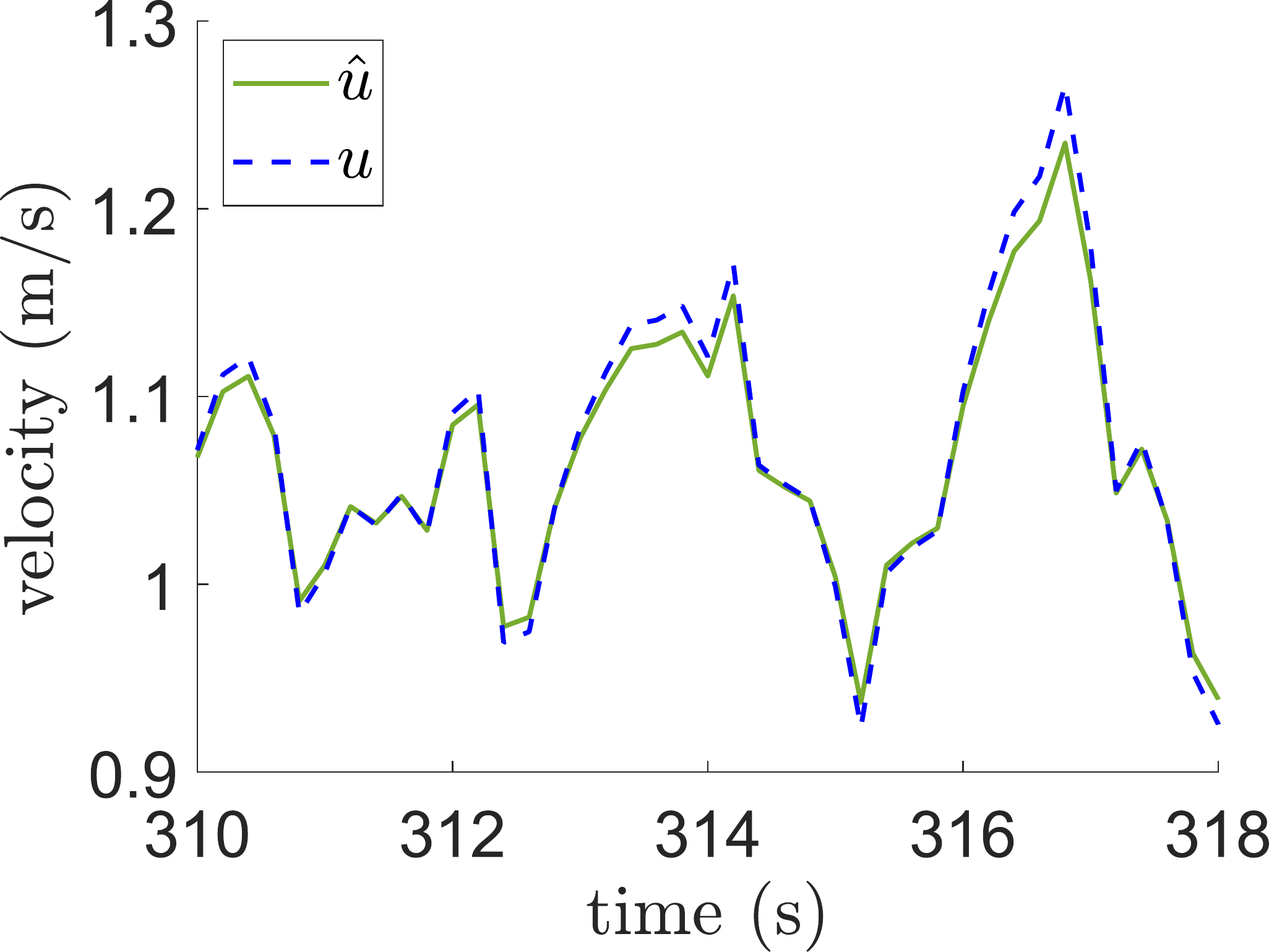}} 
    \subfigure[]{\includegraphics[width=0.32\textwidth]{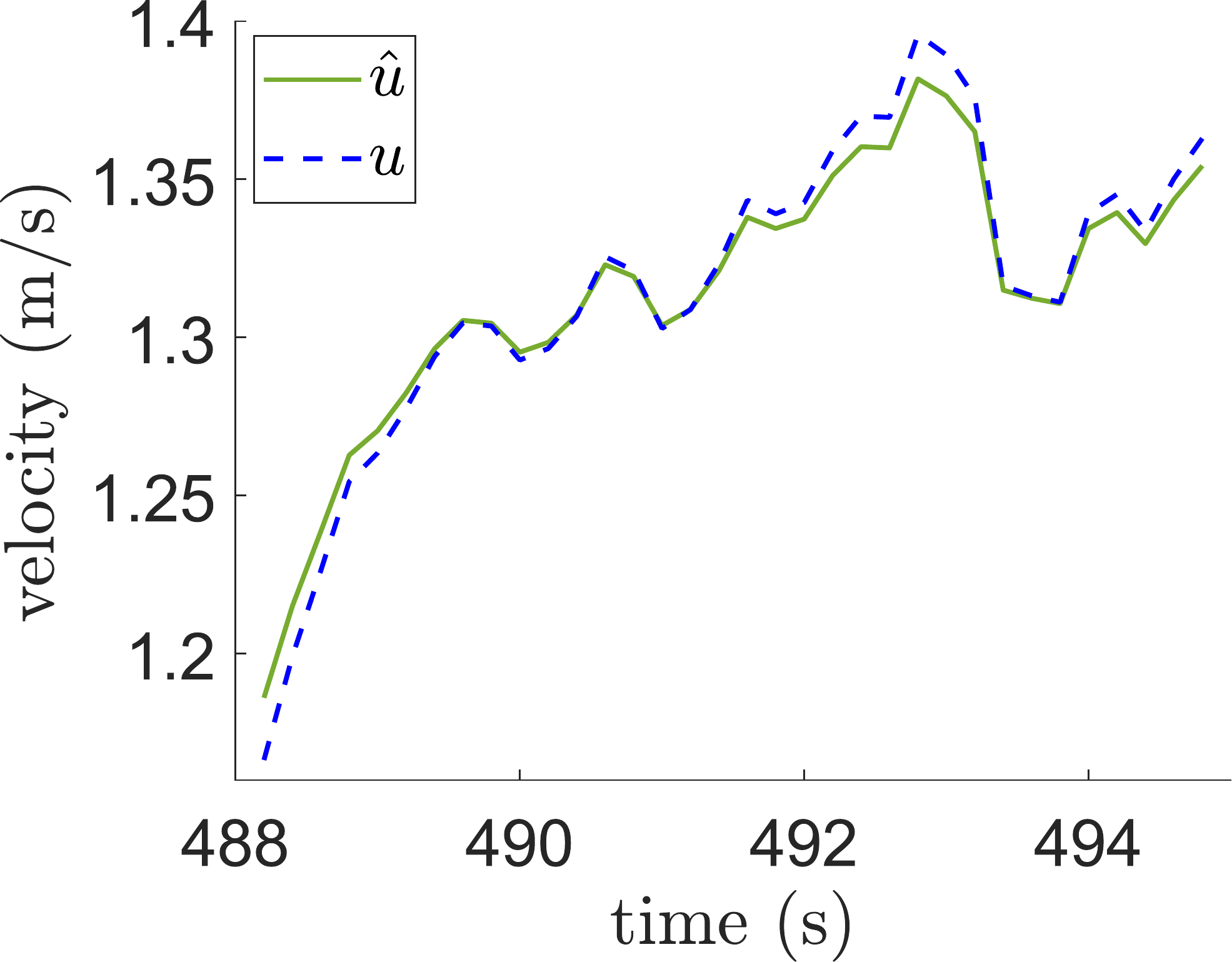}}
    \caption{Surge Velocity $u$ and Estimated Velocity $\hat{u}$ in the Validation Subset Using Second-Order Static Model}
    \label{fig_uVal}
\end{figure}

Figure~\ref{fig_uVal} shows three portions of the validation set, illustrating the estimated surge velocity during different manoeuvres, providing a wide range of velocity values, and the performance of the parameters obtained with the estimated velocity.

\begin{figure}[h!]
    \centering
    \subfigure[]{\includegraphics[width=0.32\textwidth]{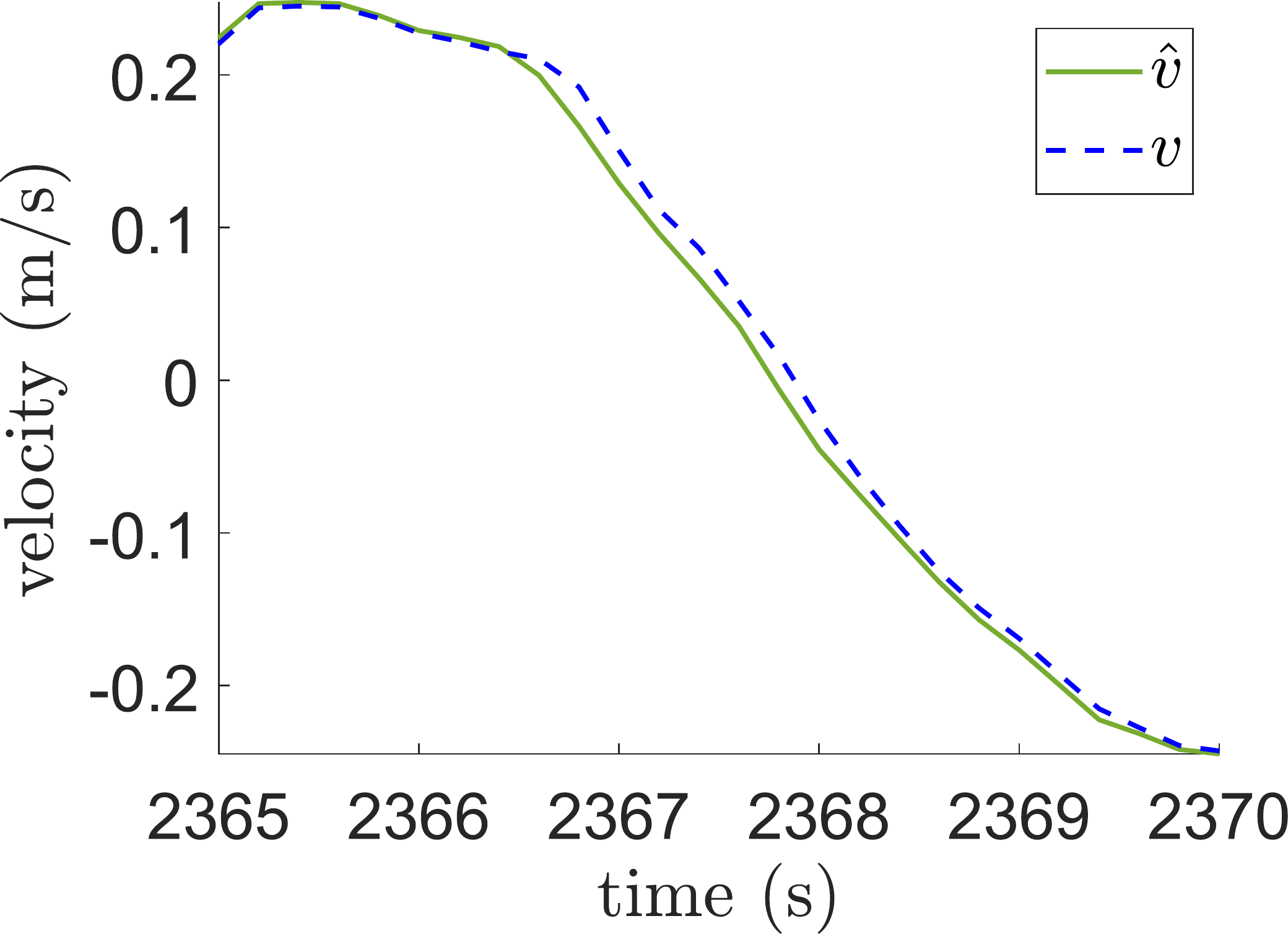}} 
    \subfigure[]{\includegraphics[width=0.32\textwidth]{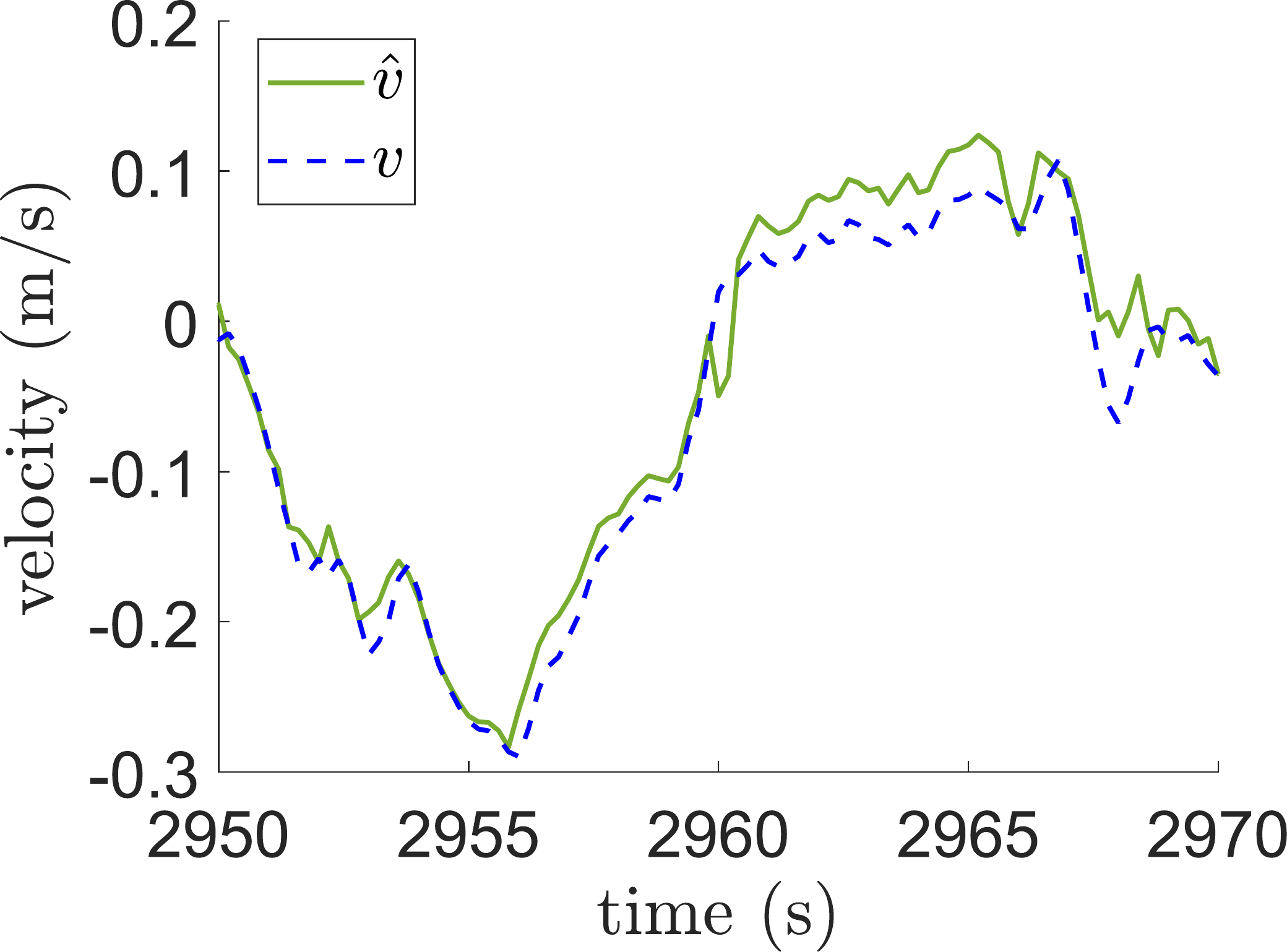}} 
    \subfigure[]{\includegraphics[width=0.32\textwidth]{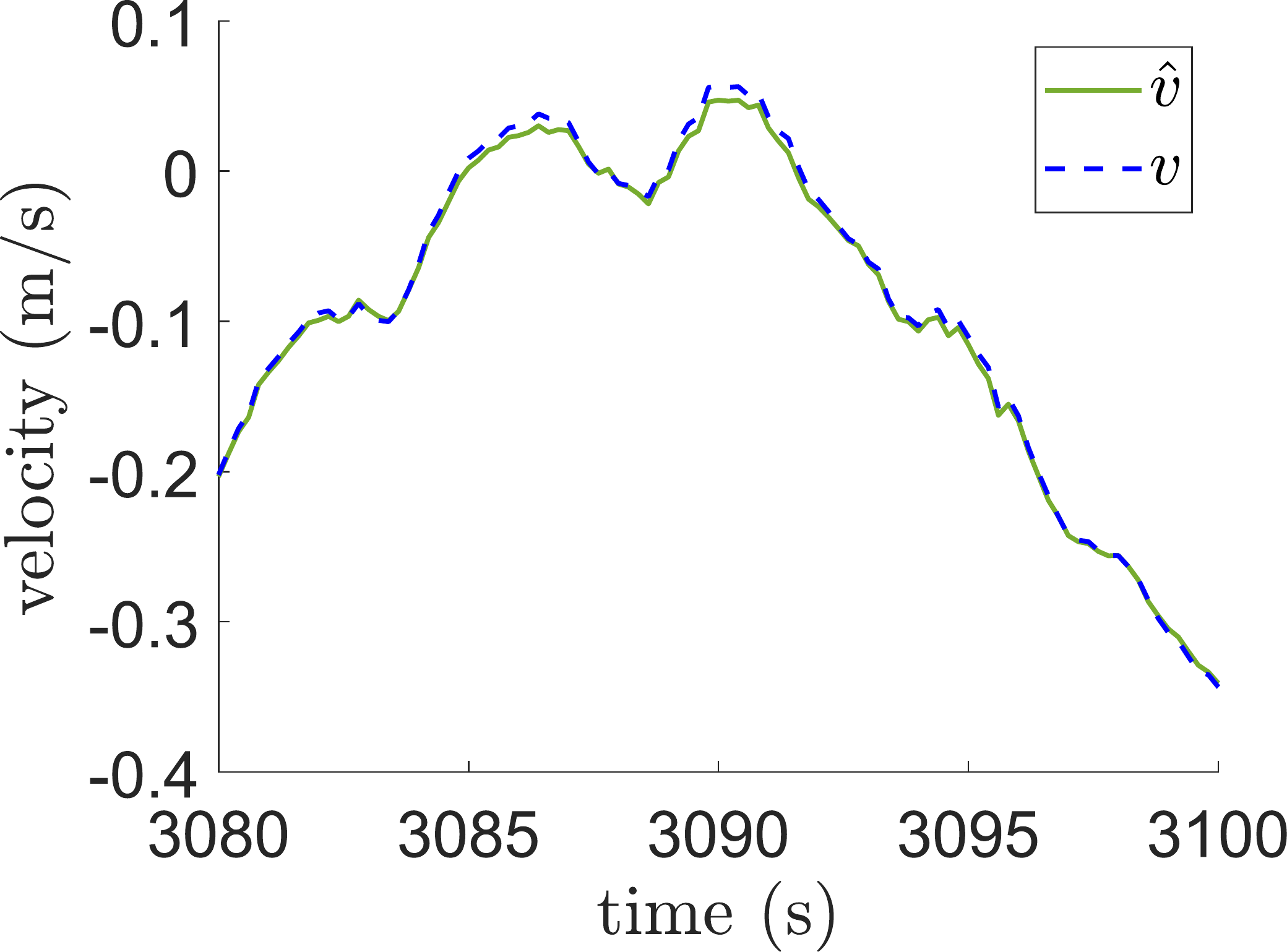}}
    \caption{Sway Velocity $v$ and Estimated Velocity $\hat{v}$ in the Validation Subset Using Second-Order Static Model}
    \label{fig_vVal}
\end{figure}

\begin{figure}[h!]
    \centering
    \subfigure[]{\includegraphics[width=0.32\textwidth]{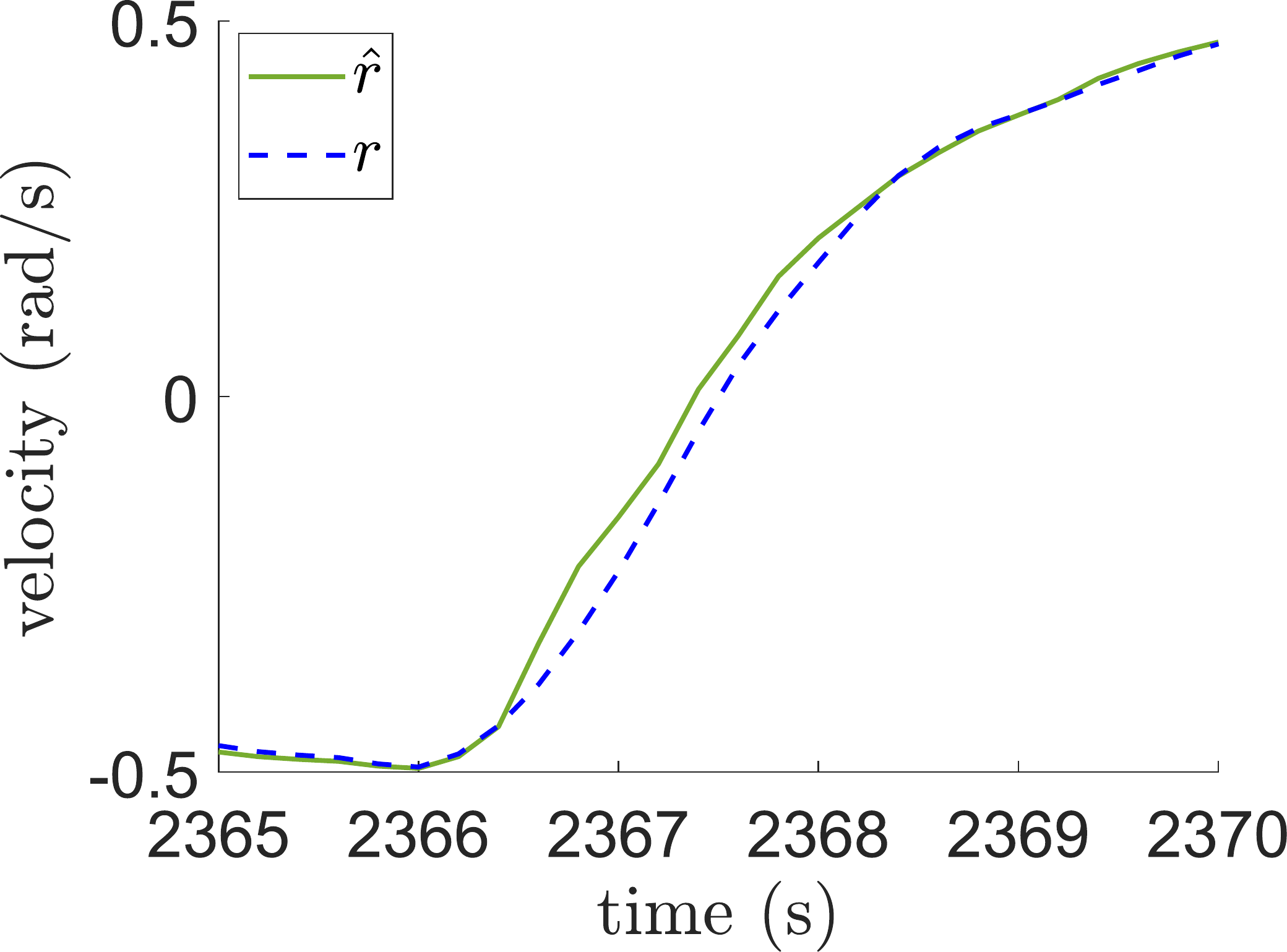}} 
    \subfigure[]{\includegraphics[width=0.32\textwidth]{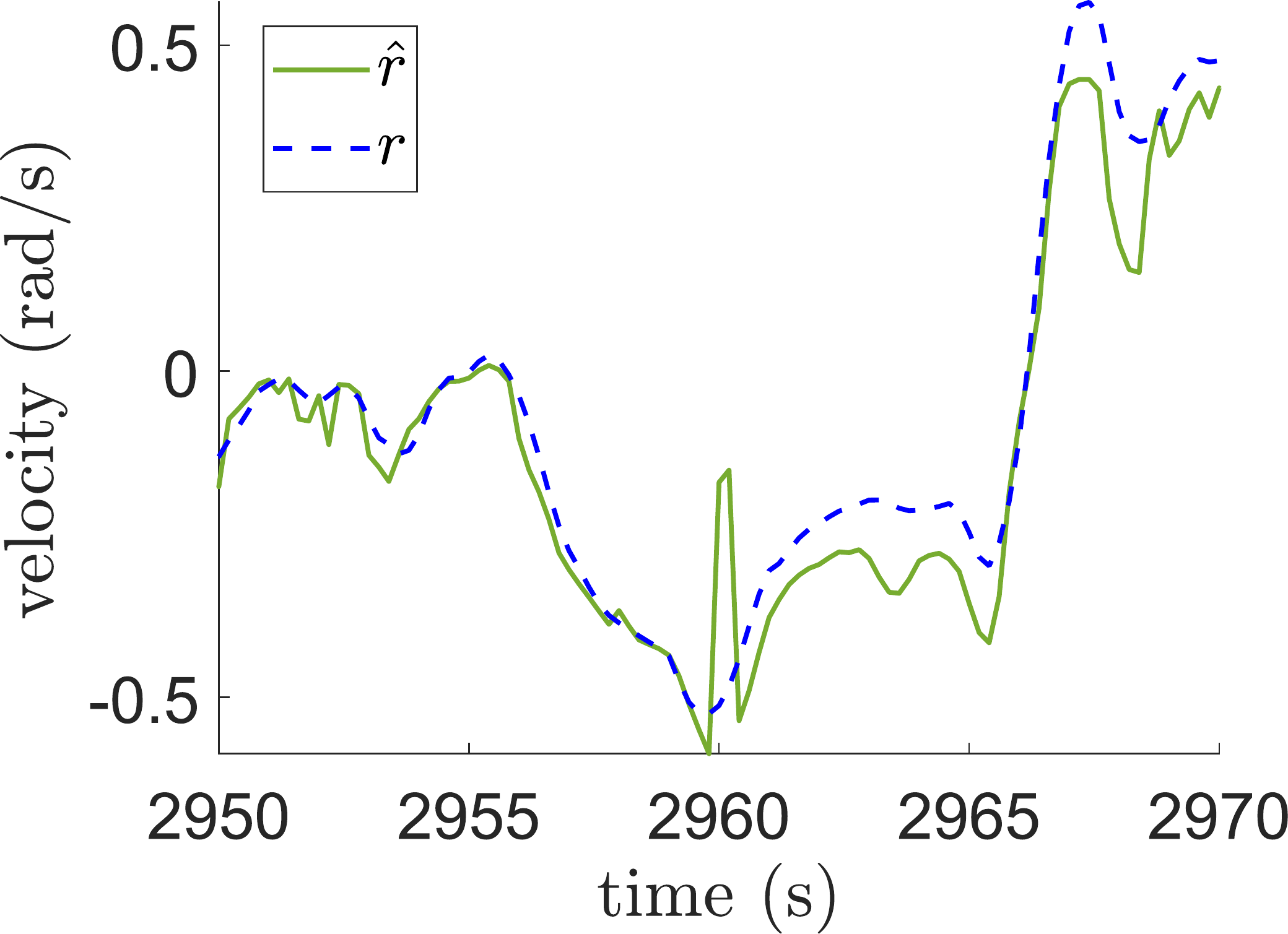}} 
    \subfigure[]{\includegraphics[width=0.32\textwidth]{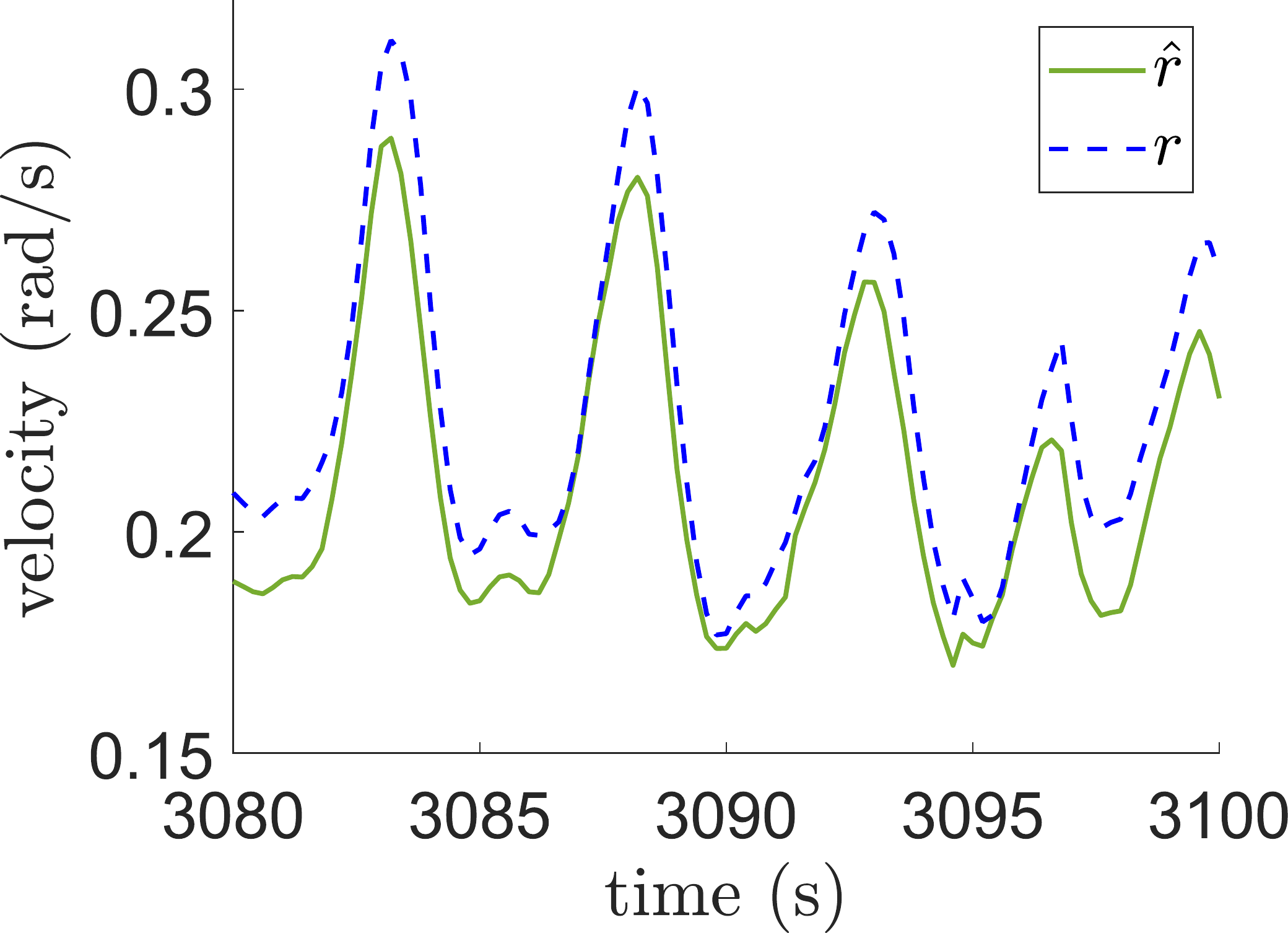}}
    \caption{Yaw Rate $r$ and Estimated Yaw Rate $\hat{r}$ in the Validation Subset Using Second-Order Static Model}
    \label{fig_rVal}
\end{figure}

Figures~\ref{fig_vVal} and~\ref{fig_rVal} present the sway velocity and the yaw rate, respectively. In the yaw rate results, the zoomed-in section of the figure reveals that the estimated values do not reach the maximum levels during sharp turns. This can still be considered a valid response, as the experimental data used for this identification involved more abrupt movements than typical operations to capture the full range of the ASV movements.

Note that the time windows shown for the graphical results represent a small percentage of the total data obtained for the validation subset. Therefore, they should be considered alongside the previously stated numerical results, which describe the errors for the entire analysed dataset.

%%%%%%%%%%%%%%%%%%%%%%%%%%%%%%%%%%%%%%%%%%%%%%%%%%%%%%%%%%%%%%%%%%%%%%%%%
\subsection{Identification Results for Dynamic Propeller Modelling}
%%%%%%%%%%%%%%%%%%%%%%%%%%%%%%%%%%%%%%%%%%%%%%%%%%%%%%%%%%%%%%%%%%%%%%%%%

The same procedure is used to initially evaluate the mean values and standard deviation of 20 randomly grouped values and is applied to the first-order dynamic model. However, in this case, since it is a dynamic model, the independent points take into account both the current values and those from the previous time step to estimate the next one. The results presented in Table~\ref{tab:MAE_R2_dyn} are quite positive and show very low standard deviation values. Consequently, for the subsequent analysis, only one group of values randomly grouped together will be examined.

\begin{table}[t]
\centering
\caption{Comparison of MAE and $R^2$ Metrics for Different Data Separation Methods Using First-Order Dynamic Model}
\label{tab:MAE_R2_dyn}
\resizebox{\columnwidth}{!}{%
\begin{tabular}{lllllll}
\toprule
                 & \textbf{MAE}$_u$  & \textbf{SD}$_u$ & \textbf{MAE}$_v$ & \textbf{SD}$_v$& \textbf{MAE}$_r$& \textbf{SD}$_r$ \\ \cline{2-7} 
\textbf{By Segments} & 0.011628 & $1.2\times 10^{-3}$     & 0.009654 & $6.3\times 10^{-4}$       &  0.007861 & $5.5 \times^{-4}$      \\
\textbf{By Independent Points}  &  \bm{$0.011770$}  & \bm{$1.1\times 10^{-4}$}     & \bm{$0.009414$} & \bm{$1.0\times 10^{-4}$}       &\bm{$0.007823$} & \bm{$6.3\times 10^{-5}$}       \\ \hline
                 & \bm{$R^2_u$}& \textbf{SD}$_u$   & \bm{$R^2_v$ }& \textbf{SD}$_v$ & \bm{$R^2_r$ } & \textbf{SD}$_r$\\ \cline{2-7} 
\textbf{By Segments} & 0.997505 & $1.5\times 10^{-3}$      & 0.990236 & $1.5\times 10^{-3}$       &   0.999523 & $1.8\times 10^{-4}$   \\
\textbf{By Independent Points}   & \bm{$0.997776$} & \bm{$6.5\times 10^{-5}$}      & \bm{$0.991055$} & \bm{$3.3\times 10^{-4}$}       &\bm{$0.999551$} & \bm{$1.0\times 10^{-5}$}    \\ \toprule
\end{tabular}%
}
\end{table}

In conducting the analysis to determine sensitivity to the number of points in the training dataset, Tables~\ref{tab:MAE_percentages_dyn} and~\ref{tab:R2_percentages_dyn} present the MAE and $R^2$ results for the various percentages of points using the dynamic model.

Both the MAE and $R^2$ values show minimal fluctuation in the different percentages, demonstrating the consistent performance of the model. Compared to the dynamic model, it is evident that using a higher percentage of data for training produces consistently better results.

\begin{table}[t]
\centering
\caption{MAE Values for Different Training Percentages Using First-Order Dynamic Model}
\label{tab:MAE_percentages_dyn}
\resizebox{\columnwidth}{!}{%
\begin{tabular}{cllllcllll}
\toprule
\multicolumn{1}{l}{\textbf{Training \%}} & \textbf{MAE}$_u$      & \textbf{MAE}$_v$      & \textbf{MAE}$_r$      &  & \multicolumn{1}{l}{\textbf{Validation \%}} & \textbf{MAE}$_u$      & \textbf{MAE}$_v$      & \textbf{MAE}$_r$      &  \\ \cline{1-4} \cline{6-9}
70                              & \bm{$0.011637$} & \bm{$0.010785$} & \bm{$0.008415$} &  & 30                                & \bm{$0.033513$} & \bm{$0.007386$} & \bm{$0.016770$} &  \\
60                              & 0.012871 & 0.011076 & 0.008876 &  & 30                                & 0.034702 & 0.008325 & 0.030219 &  \\
50                              & 0.015348 & 0.012492 & 0.008611 &  & 30                                & 0.043901 & 0.008154 & 0.022668 &  \\
\bottomrule
\end{tabular}
}
\end{table}

\begin{table}[h]
\centering
\caption{$R^2$ Values for Different Training Percentages Using First-Order Dynamic Model}
\label{tab:R2_percentages_dyn}
\resizebox{\columnwidth}{!}{%
\begin{tabular}{cllllcllll}
\toprule
\textbf{Training \%} & \bm{$R^2_u$ }     & \bm{$R^2_v$ }      & \bm{$R^2_r$}       &  & \textbf{Validation \%} & \bm{$R^2_u$ }      & \bm{$R^2_v$}       & \bm{$R^2_r$}       &  \\ \cline{1-4} \cline{6-9}
70       & \bm{$0.998445$} & 0.987943 & \bm{$0.999573$} &  & 30         & \bm{$0.943559$} & \bm{$0.996225$} & \bm{$0.997599$} &  \\
60       & 0.996819 &\bm{$0.988532$} & 0.999137 &  & 30         & 0.939760 & 0.995190 & 0.992966 &  \\
50       & 0.989830 & 0.985986 & 0.999189 &  & 30         & 0.909676 & 0.995331 & 0.996110 &  \\ \toprule
\end{tabular}%
}
\end{table}

Observing the graphical analysis, the same time windows are used with identical training and validation datasets as for the previous model analysed. This approach enables a thorough comparison and evaluation of both models analysed in this section. 

\begin{figure}[h]
    \centering
    \subfigure[]{\includegraphics[width=0.32\textwidth]{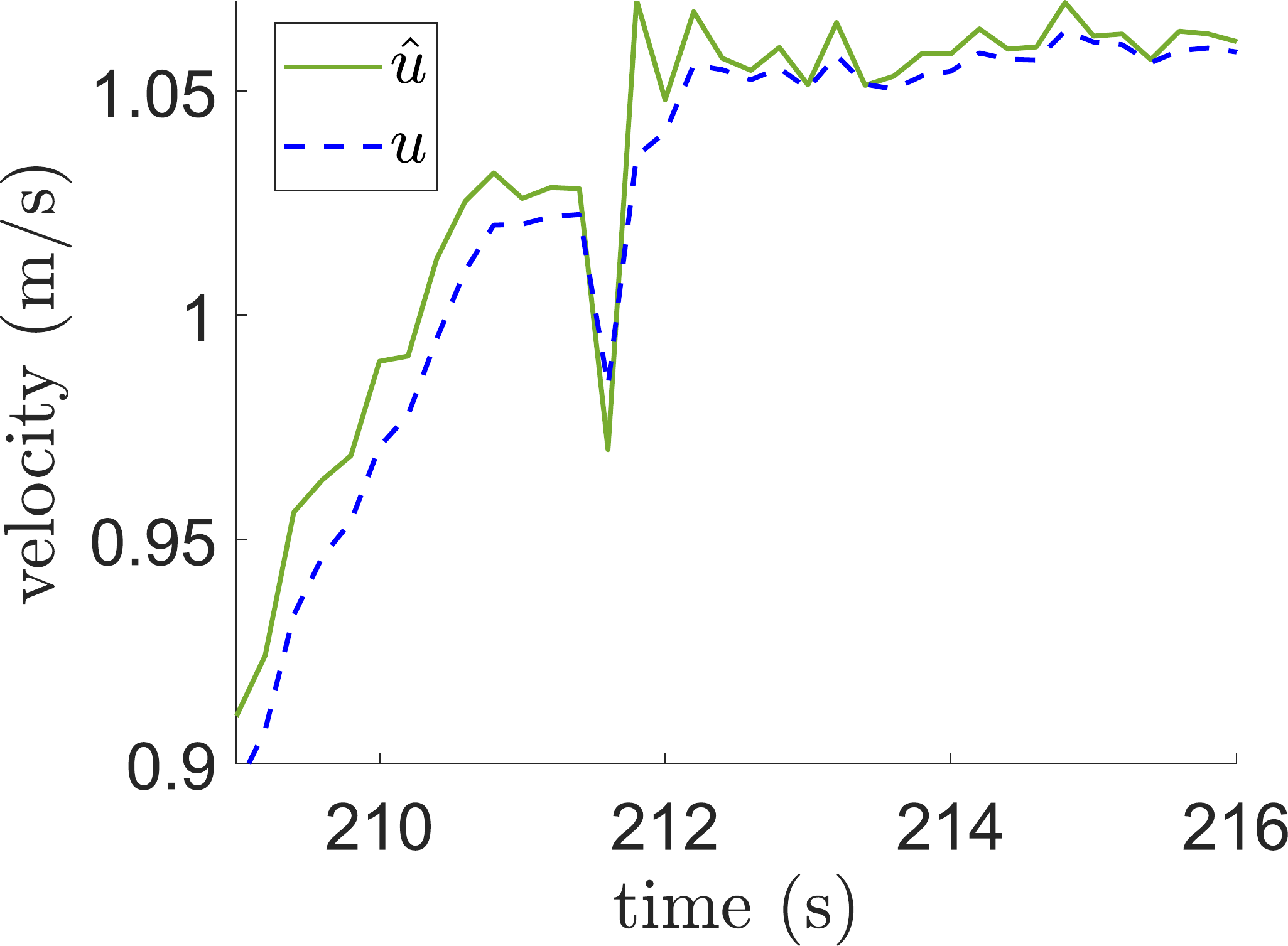}} 
    \subfigure[]{\includegraphics[width=0.32\textwidth]{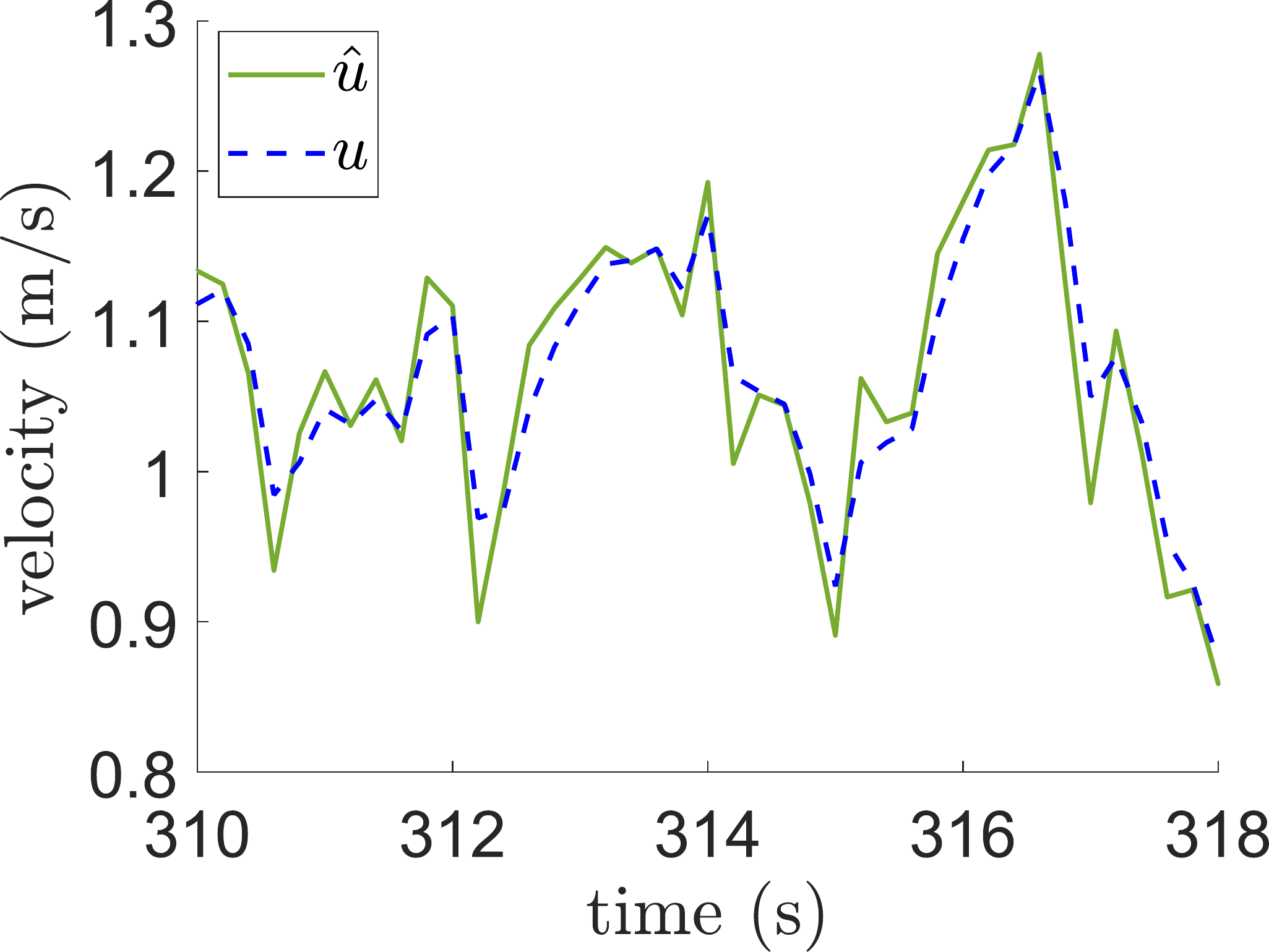}} 
    \subfigure[]{\includegraphics[width=0.32\textwidth]{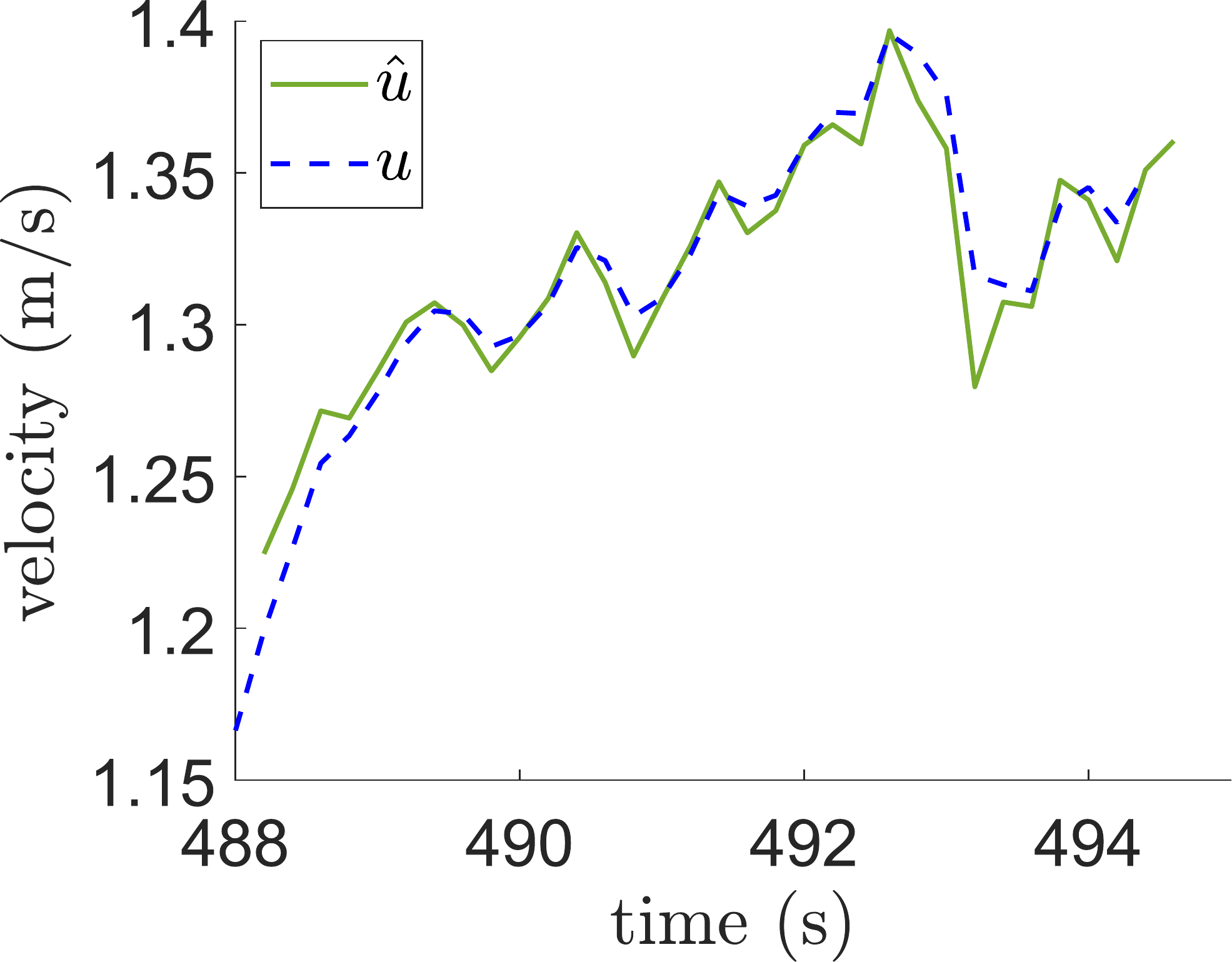}}
    \caption{Surge Velocity $u$ and Estimated Velocity $\hat{u}$ in the Validation Subset Using First-Order Dynamic Model}
    \label{fig_uVal_dyn}
\end{figure}

\begin{figure}[h]
    \centering
    \subfigure[]{\includegraphics[width=0.32\textwidth]{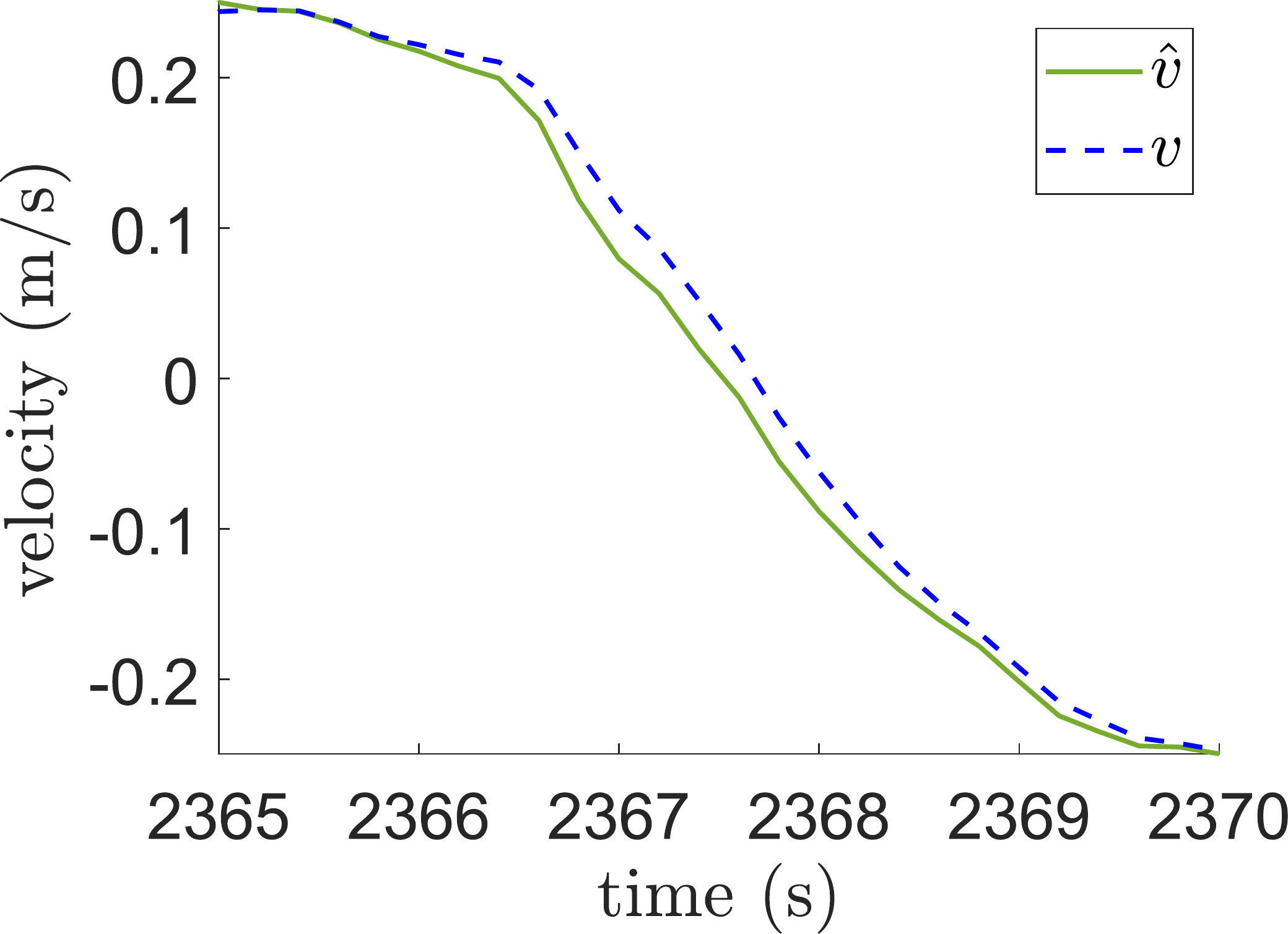}} 
    \subfigure[]{\includegraphics[width=0.32\textwidth]{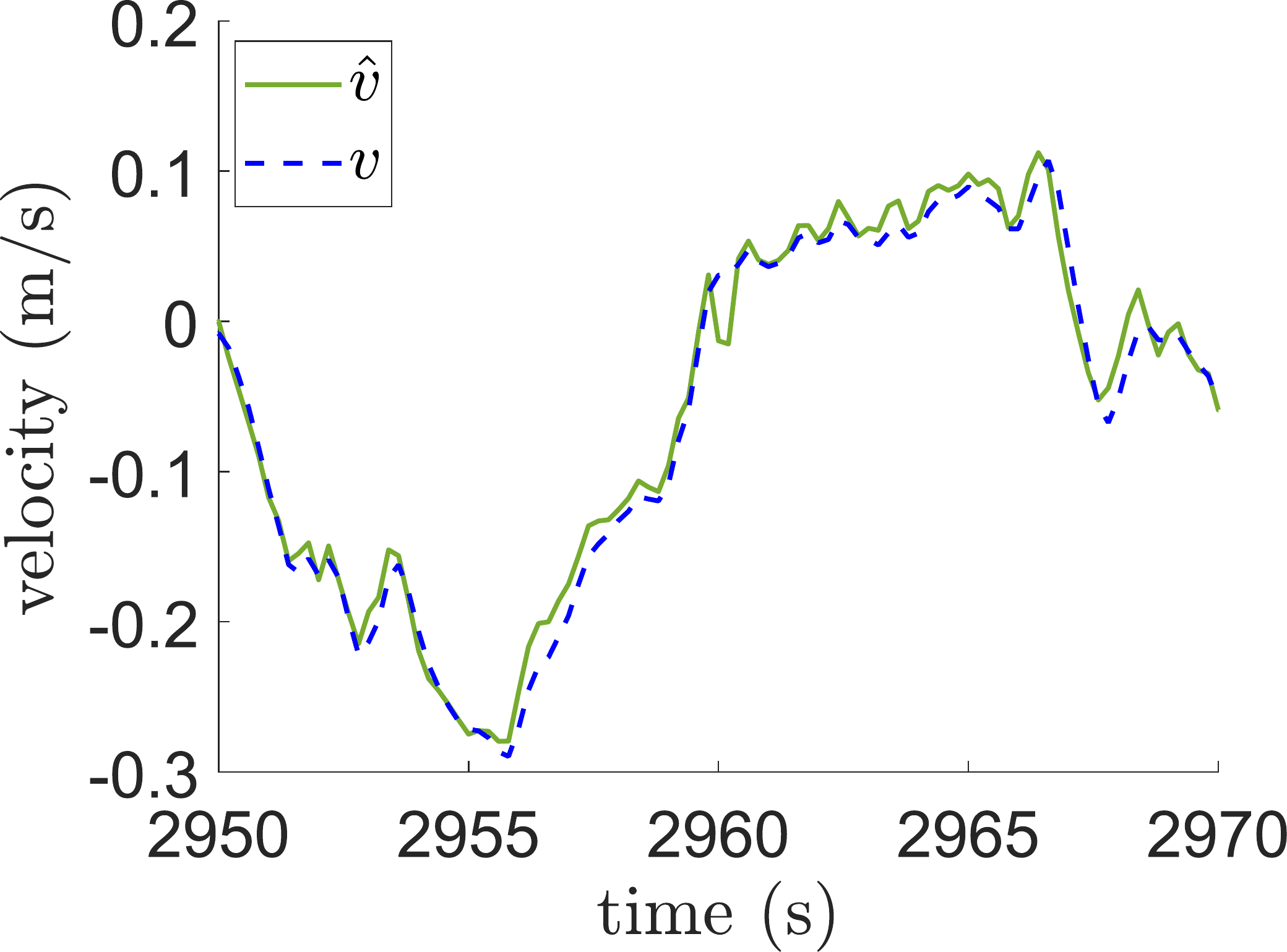}} 
    \subfigure[]{\includegraphics[width=0.32\textwidth]{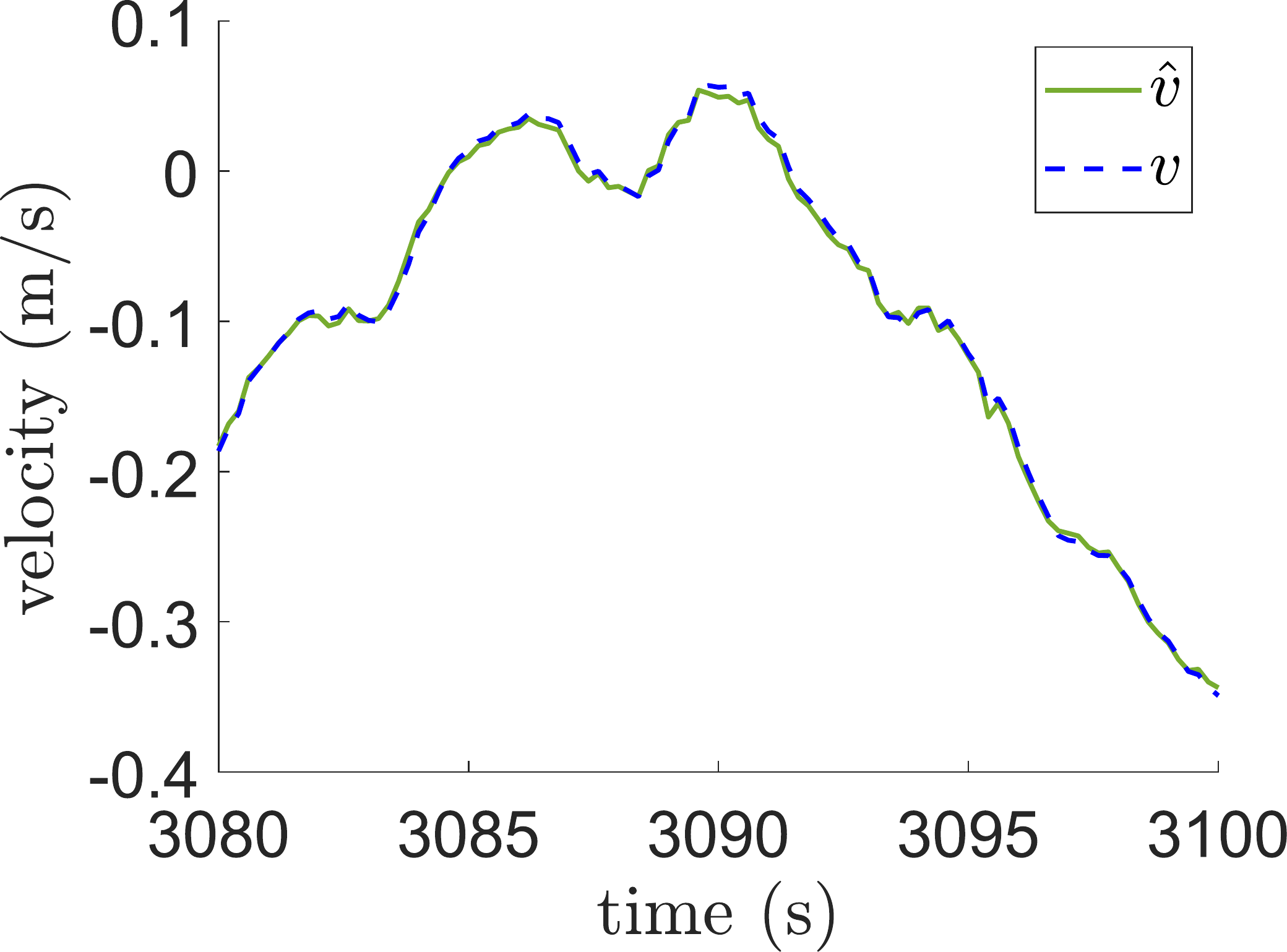}}
    \caption{Sway Velocity $v$ and Estimated Velocity $\hat{v}$ in the Validation Subset Using First-Order Dynamic Model}
    \label{fig_vVal_dyn}
\end{figure}

\begin{figure}[h!]
    \centering
    \subfigure[]{\includegraphics[width=0.32\textwidth]{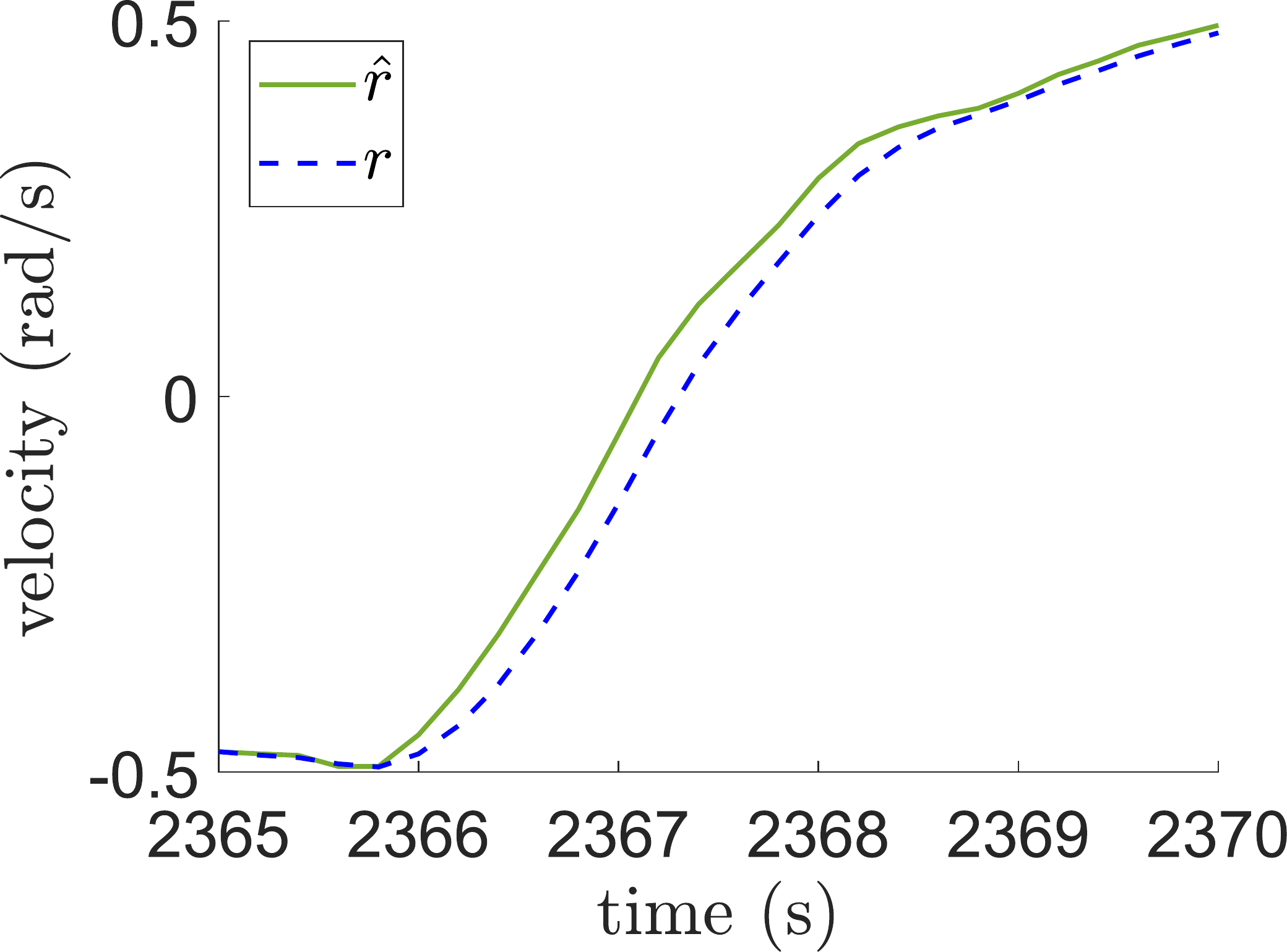}} 
    \subfigure[]{\includegraphics[width=0.32\textwidth]{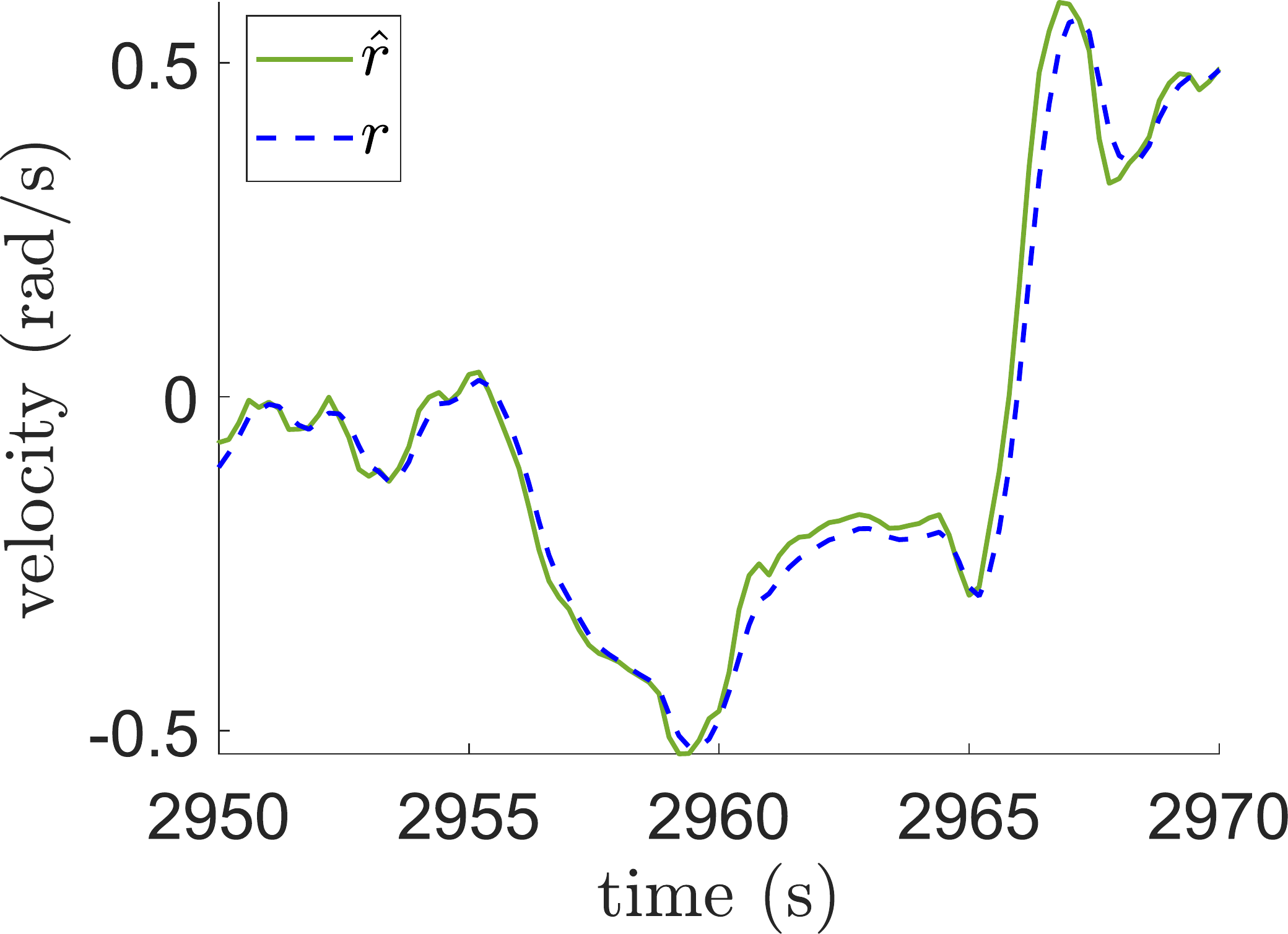}} 
    \subfigure[]{\includegraphics[width=0.32\textwidth]{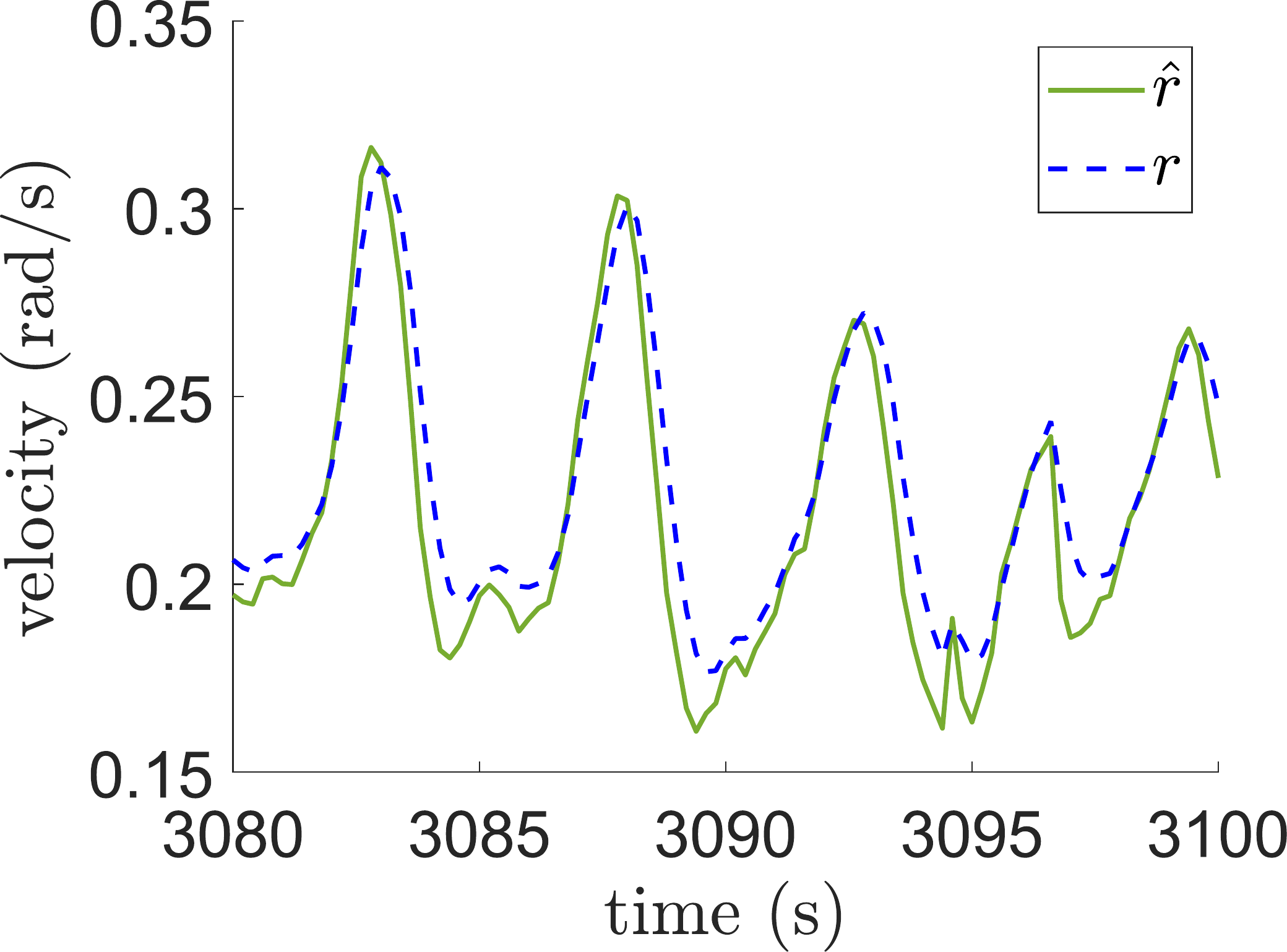}}
    \caption{Yaw Rate $r$ and Estimated Yaw Rate $\hat{r}$ in the Validation Subset Using First-Order Dynamic Model}
    \label{fig_rVal_dyn}
\end{figure}

Figure~\ref{fig_uVal_dyn} displays the estimated surge velocity, and Figures~\ref{fig_vVal_dyn} and~\ref{fig_rVal_dyn} show the results for sway velocity and yaw rate, respectively. 

Regarding the surge results, the estimated values exhibit pronounced peaks and appear less smooth, as reflected in the metrics previously presented.

Considering sway velocity and yaw rate, the model demonstrates good and smooth tracking of the estimated velocities, although some peaks are more pronounced than those of the actual velocities, particularly in regions with abrupt velocity changes.

To obtain the estimated velocity values, it is necessary to determine the value of $\alpha$ using the identification equations for the three velocities~\eqref{eq_alpha}. The resulting value of $\alpha$ is 0.998, which ensures the stability of the model, as discussed in Remark~\ref{remark:alpha}.

%%%%%%%%%%%%%%%%%%%%%%%%%%%%%%%%%%%%%%%%%%%%%%%%%%%%%%%%%%%%%%%%%%%%%
\subsection{Discussion and Comparison of the Static and Dynamic Models}
%%%%%%%%%%%%%%%%%%%%%%%%%%%%%%%%%%%%%%%%%%%%%%%%%%%%%%%%%%%%%%%%%%%%%

Finally, the differences between static and dynamic models are discussed. In particular, the first-order dynamic model produces superior results for the sway velocity and yaw rate, according to the metrics employed. Furthermore, graphical analysis shows improved accuracy in the estimated yaw rate, particularly in sharper curves, although some minor peaks persist compared to the second-order model.

In contrast, the surge velocity shows better results when using the static model. To understand this difference, it is important to note that the dataset used to identify the surge dynamics parameters is smaller than that for the sway and yaw. Furthermore, in the dynamic model, a single value of $\alpha$ is applied at all three velocities, which seems to favour the dynamics of the sway and yaw.

The static and dynamic models perform well, with the dynamic model offering slightly better results; however, it is more complex to implement and requires the stability of the dynamic model of the propellers (with $\alpha < 1$). Therefore, choosing which model to use ultimately depends on the preference and requirements of the user.

To conclude, Tables~\ref{tab:static_values} and \ref{tab:dynamic_values} present the values of the parameters identified for the static second-order model, as well as the values obtained using the dynamic first-order model. 

\begin{table}[h!]
\centering
\caption{Identified Parameters and Corresponding Units of Measurement Used for the Second-Order Static Model}
\label{tab:static_values}
\resizebox{0.6\columnwidth}{!}{%
\begin{tabular}{ccccccc}
\toprule
      & \multicolumn{6}{c}{\textbf{Second-Order Static Model}}                                                                                         \\ \hline
Index & $X_u$                   & Unit                      & $X_v$      & Unit                      & $X_r$      & Unit                      \\ \hline
1     & -0.0616              & (m/s)$^-1$ & -0.0214 & (m/s)$^-1$ & 0.0424  & (m/s)$^-1$ \\
2     & 0.0272               & (m/s)$^-1$ & -0.0714 & (m/s)$^-1$ & -0.1168 & (m/s)$^-1$ \\
3     & -0.0098              & (m/s)$^-1$ & -0.0334 & (m/s)$^-1$ & 0.0369  & (m/s)$^-1$ \\
4     & 0.0072               & -                         & -0.0123 & (m/s)$^-1$ & -0.0077 & (m/s)$^-1$ \\
5     & -0.0001              & m/s                       & -0.0400 & (m/s)$^-1$ & -0.0310 & (m/s)$^-1$ \\
6     & \bm{$-0.0145$}              & m/s                       & 0.0329  & (m/s)$^-1$ & -0.0931 & (m/s)$^-1$ \\
7     & \bm{$0.1403$}               & m/s                       & 0.0407  & -                         & 0.0882  & -                         \\
8     & \multicolumn{1}{l}{} & \multicolumn{1}{l}{}      & 0.0334  & -                         & -0.1505 & -                         \\
9     & \multicolumn{1}{l}{} & \multicolumn{1}{l}{}      & -0.0018 & m/s                       & 0.0055  & m/s                       \\
10    & \multicolumn{1}{l}{} & \multicolumn{1}{l}{}      & \bm{$0.0092$}  & m/s                       & \bm{$-0.0872$} & m/s                       \\
11    & \multicolumn{1}{l}{} & \multicolumn{1}{l}{}      & \bm{$-0.0381$} & m/s                       & \bm{$-0.0743$} & m/s                       \\
12    & \multicolumn{1}{l}{} & \multicolumn{1}{l}{}      & \bm{$-0.0032$} & m/s                       & \bm{$0.0814$}  & m/s                       \\
13    & \multicolumn{1}{l}{} & \multicolumn{1}{l}{}      & \bm{$-0.0505$} & m/s                       & \bm{$0.3206$}  & m/s                       \\ \bottomrule
\end{tabular}%
}
\end{table}

The tables highlight the values required for the input gain equations presented in~\eqref{eq_IG:u} and~\eqref{eq:IG_vr}, using the vessel \emph{Yellowfish}, where the dynamic model also incorporates the previously mentioned $\alpha$ value. All identified parameter values are presented to enable replication of the results.

\begin{table}[h!]
\centering
\caption{Identified Parameters and Corresponding Units of Measurement Used for the First-Order Dynamic Model}
\label{tab:dynamic_values}
\resizebox{0.6\columnwidth}{!}{%
\begin{tabular}{ccccccc}
\toprule
\multicolumn{7}{c}{ \textbf{First-Order Dynamic Model Parameters}}                                                                                   \\ \hline
Index & $X_{u,dyn}$   & Unit                      & $X_{v,dyn}$   & Unit                      & $X_{v,dyn}$   & Unit                      \\ \hline
1     & 0.8062  & -                         & 0.6904  & -                         & 0.9886  & -                         \\
2     & 0.0893  & (m/s)$^-1$ & -0.0091 & (m/s)$^-1$ & 0.0577  & (m/s)$^-1$ \\
3     & -0.1095 & (m/s)$^-1$ & 0.0237  & (m/s)$^-1$ & 0.1845  & (m/s)$^-1$ \\
4     & -0.0341 & (m/s)$^-1$ & -0.0261 & (m/s)$^-1$ & -0.0190 & (m/s)$^-1$ \\
5     & -0.7725 & -                         & -0.0005 & (m/s)$^-1$ & -0.0054 & (m/s)$^-1$ \\
6     & -0.1395 & (m/s)$^-1$ & 0.4328  & (m/s)$^-1$ & 0.0667  & (m/s)$^-1$ \\
7     & 0.1622  & (m/s)$^-1$ & 0.3075  & (m/s)$^-1$ & -0.0020 & (m/s)$^-1$ \\
8     & 0.0301  & (m/s)$^-1$ & -0.6978 & -                         & -0.2723 & -                         \\
9     & 0.0016  & m/s                       & -0.1285 & -                         & -1.0430 & -                         \\
10    & \bm{$-0.0106$} & m/s                       & -0.0085 & (m/s)$^-1$ & -0.1305 & (m/s)$^-1$ \\
11    & \bm{$0.0611$}  & m/s                       & -0.0494 & (m/s)$^-1$ & -0.2120 & (m/s)$^-1$ \\
12    &         &                           & 0.0081  & (m/s)$^-1$ & -0.0059 & (m/s)$^-1$ \\
13    &         &                           & 0.0004  & (m/s)$^-1$ & 0.0577  & (m/s)$^-1$ \\
14    &         &                           & -0.4353 & (m/s)$^-1$ & -0.1116 & (m/s)$^-1$ \\
15    &         &                           & -0.3078 & (m/s)$^-1$ & 0.0358  & (m/s)$^-1$ \\
16    &         &                           & 0.1478  & -                         & 0.3419  & -                         \\
17    &         &                           & -0.0007 & m/s                       & -0.0013 & m/s                       \\
18    &         &                           & \bm{$-0.007$}  & m/s                       & \bm{$-0.0422$} & m/s                       \\
19    &         & \multicolumn{1}{l}{}      & \bm{$-0.007$}  & m/s                       & \bm{$-0.0385$} & m/s                       \\
20    &         & \multicolumn{1}{l}{}      & \bm{$0.0101$}  & m/s                       & \bm{$0.0394$}  & m/s                       \\
21    &         & \multicolumn{1}{l}{}      & \bm{$-0.0214$} & m/s                       & \bm{$0.0168$}  & m/s                       \\ \bottomrule
\end{tabular}%
}
\end{table}

%%%%%%%%%%%%%%%%%%%%%%%%%%%%%%%%%%%%%%%%%%%%%%%%%
\section{Conclusions}
\label{sec:sample:Conclusions}
%%%%%%%%%%%%%%%%%%%%%%%%%%%%%%%%%%%%%%%%%%%%%%%%%

This work presents a parameter identification methodology to determine the input gain of an ASV with fully unknown coefficients. The unknown model parameters can be estimated through position and heading measurements, thus eliminating the need for acceleration data. 

This study presents a practical method for jointly identifying the propeller model and the
inertia matrix of an ASV, which means the effect of the propellers on its motion. The identification method included, with experimental results, both a second-order static model and a first-order dynamic model for the actuation system.

The models were applied to obtain the parameters of the \emph{Yellowfish} ASV, with significant emphasis on experimental data acquisition to ensure robust parameter functionality across a significant operational range. 

The effectiveness of the proposed methods is validated by applying identification metrics to the results, also allowing comparison between different models. For the second-order static model, a determination coefficient greater than 0.98 was achieved for surge velocity, sway velocity, and yaw rate. In contrast, when applying the first-order dynamic model to the same experimental data, the results yielded a range of 0.94 to 0.99 across all three velocities. These consistent values highlight the reliability and effectiveness of the models in accurately capturing vessel dynamics.

Furthermore, the importance of replicating the results is demonstrated, as it contributes to a deeper understanding of the models. To support this, the identified values obtained are reported in Tables~\ref{tab:static_values} and~\ref{tab:dynamic_values}, and all experimental data related to this work are available from~\cite{Morel2024}.

\section*{Acknowledgements}
This work was supported by the Agencia Estatal de Investigación (AEI) through the Project AQUATRONIC under Grant PID2021 126921OA-C22, and through the Project ECOPORT under Grant TED2021-131326A-C22. We also thank the Consejería de Fomento, Infraestructuras y Ordenación del Territorio of the Junta de Andalucía for their collaboration and support in conducting tests in the Parque del Alamillo.

\appendix

%%%%%%%%%%%%%%%%%%%%%%%%%%%%%%%%%%%%%%%%%%%%%%%%%%%%%%%%%%%%%%%%%
\section{Identification of the Input-Gain for the Sway and Yaw}
\label{apendix_static_swayyaw}
%%%%%%%%%%%%%%%%%%%%%%%%%%%%%%%%%%%%%%%%%%%%%%%%%%%%%%%%%%%%%%%%%
In the analysis of sway and yaw, the three regions under evaluation $(f,f)$, $(f,r)$ and $(r,f)$ are employed, using~\eqref{eq_tau_r_2} to define the torque $\tau_r$. Consequently, the sway and yaw dynamics, considering the quadratic static model for the propellers, are as follows:
\begin{align*} 
    v(k+1) &= v(k)+h\bm{M}^{-1}(2,3)\frac{d}{2}\Bigg((a_{j}^{L}-a_{j}^{R})\left(\bar{\delta}(k)^2+\frac{\Delta\delta(k)^2}{4}\right)\nonumber\\
         &\quad+(a_{j}^{L}+a_{j}^{R})\bar{\delta}(k)\Delta\delta(k) + (b_{j}^{L}-b_{j}^{R})\bar{\delta}(k)+(b_{j}^{L}+b_{j}^{R})\frac{\Delta\delta(k)}{2}\Bigg)\nonumber\\
         &\quad+ hP_{v}^{(2,2)}v(k)|v(k)|+hP_{v}^{(2,3)}v(k)|r(k)|+hP_{v}^{(3,2)}r(k)|v(k)|\nonumber\\ 
         &\quad+hP_{v}^{(3,3)}r(k)|r(k)| + 2hQ_{v}^{(1,2)}u(k)v(k)+2hQ_{u}^{(1,3)}u(k)r(k)\nonumber\\
    &\quad+hR_{v}^{(2)}v(k)+hR_{v}^{(3)}r(k)+hc_{v}(k), 
\end{align*} %\label{eq_vk1_diff} 
\begin{align*} 
        r(k+1) &=r(k)+h\bm{M}^{-1}(3,3)\frac{d}{2}\Bigg((a_{j}^{L}-a_{j}^{R})\left(\bar{\delta}(k)^2+\frac{\Delta\delta(k)^2}{4}\right)\nonumber\\
         &\quad+(a_{j}^{L}+a_{j}^{R})\bar{\delta}(k)\Delta\delta(k) + (b_{j}^{L}-b_{j}^{R})\bar{\delta}(k)+(b_{j}^{L}+b_{j}^{R})\frac{\Delta\delta(k)}{2}\Bigg)\nonumber\\
         &\quad + hP_{r}^{(2,2)}v(k)|v(k)| +hP_{r}^{(2,3)}v(k)|r(k)|+hP_{r}^{(3,2)}r(k)|v(k)|\nonumber \\
    &\quad+hP_{r}^{(3,3)}r(k)|r(k)| + 2hQ_{r}^{(1,2)}u(k)v(k)+2hQ_{r}^{(1,3)}u(k)r(k)\nonumber\\
    &\quad+hR_{r}^{(2)}v(k)+hR_{r}^{(3)}r(k)+hc_{r}(k), %\label{eq_rk1_diff}
\end{align*}

where $a_{j}^{i},b_{j}^{i}$ with $i\in\{L,R\},j\in\{f,r\}$ are the forward and reverse parameters to be identified for the propellers.

The unknown coefficients appear linearly, easing the identification. The regressor equations for sway in matrix form, this is, $A_vX_v=b_v$, can be derived with:
\begin{align*}
    X_v &= h\Bigg[P_{v}^{(2,2)}~~P_{v}^{(2,3)}~~P_{v}^{(3,2)}~~P_{v}^{(3,3)} ~~2Q_v^{(1,2)}~~2Q_v^{(1,3)}~~R_v^{(2)}~~R_v^{(3)}~~\Bar{c}_v \nonumber \\ 
    &\quad~~\bm{M}^{-1}(2,3)\frac{d}{2}(a_{j}^{L}-a_{j}^{R})~~\bm{M}^{-1}(2,3)\frac{d}{2}(a_{j}^{L}+a_{j}^{R})\nonumber \\
    &\quad~~\bm{M}^{-1}(2,3)\frac{d}{2}(b_{j}^{L}-b_{j}^{R})~~\bm{M}^{-1}(2,3)\frac{d}{2}(b_{j}^{L}+b_{j}^{R})\Bigg]^T, \\ 
    A_v(k) &= \Bigg[v(k)|v(k)|~~v(k)|r(k)|~~r(k)|v(k)|~~~r(k)|r(k)|\nonumber \\
    &\quad~~u(k)v(k)~~u(k)r(k)~~v(k)~~r(k)~~1~~\nonumber \\
    &\quad \left(\bar{\delta}(k)^2+\frac{\Delta\delta(k)^2}{4}\right)~~\bar{\delta}(k)\Delta\delta(k)~~\bar{\delta}(k)~~\frac{\Delta\delta(k)}{2}\Bigg], \\
    b_v(k) &= v(k+1) - v(k).
\end{align*}

The same process is carried out for the yaw, resulting in the following equations:
\begin{align*}
    X_r &= h\Bigg[P_{r}^{(2,2)}~~P_{r}^{(2,3)}~~P_{r}^{(3,2)}~~P_{r}^{(3,3)} ~~2Q_r^{(1,2)}~~2Q_r^{(1,3)}~~R_r^{(2)}~~R_r^{(3)}~~\Bar{c}_r \nonumber \\ 
    &\quad~~\bm{M}^{-1}(3,3)\frac{d}{2}(a_{j}^{L}-a_{j}^{R})~~\bm{M}^{-1}(3,3)\frac{d}{2}(a_{j}^{L}+a_{j}^{R})\nonumber \\
    &\quad~~\bm{M}^{-1}(3,3)\frac{d}{2}(b_{j}^{L}-b_{j}^{R})~~\bm{M}^{-1}(3,3)\frac{d}{2}(b_{j}^{L}+b_{j}^{R})\Bigg]^T, \\
    A_r(k) &= \Bigg[v(k)|v(k)|~~v(k)|r(k)|~~r(k)|v(k)|~~~r(k)|r(k)|\nonumber \\ &\quad~~u(k)v(k)~~u(k)r(k)~~v(k)~~r(k)~~1~~\nonumber \\
    &\quad \left(\bar{\delta}(k)^2+\frac{\Delta\delta(k)^2}{4}\right)~~\bar{\delta}(k)\Delta\delta(k)~~\bar{\delta}(k)~~\frac{\Delta\delta(k)}{2}\Bigg], \\
    b_r(k) &= r(k+1) - r(k).
\end{align*}

It is important to note that in the $(f,f)$ case, all four terms $X_v(10)$, $X_v(12)$ for the sway, and $X_r(10)$, $X_r(12)$ for the yaw are cancelled out because the left and right propellers have the same coefficient.

%%%%%%%%%%%%%%%%%%%%%%%%%%%%%%%%%%%%%%%%%%%%%%%%%%%%%%%%%%%%%%%%%%%%%%%%%%%%
 \section{Input-Gain Identification for Dynamic Propeller Modelling}
 \label{apendix_dynamic}
 %%%%%%%%%%%%%%%%%%%%%%%%%%%%%%%%%%%%%%%%%%%%%%%%%%%%%%%%%%%%%%%%%%%%%%%%%%%%
 Consider a dynamic propeller model, with the left and right thrusts $T(\delta(k))$ being described as~\eqref{eq_prop_dyn_model}. Then, according to~\eqref{eq_F_u_tau_r}, the force $F_u$ and torque $\tau_r$ can be written as follows:
\begin{align*}
    F_u(\bar{\delta}(k),\Delta\delta(k)) &= \alpha \left( T(\delta_L(k-1) + T(\delta_R(k-1) )\right) + \nonumber\\
       &\quad \beta \left((a_{j}^{L}+a_{j}^{R})\left(\bar{\delta}(k-1)^2+\frac{\Delta\delta(k-1)^2}{4}\right) \right. \nonumber\\
       &\quad \left.+(a_{j}^{L}-a_{j}^{R})\bar{\delta}(k-1)\Delta\delta(k-1) \right. \nonumber\\
       &\quad \left. + (b_{j}^{L}+b_{j}^{R})\bar{\delta}(k-1)+(b_{j}^{L}-b_{j}^{R})\frac{\Delta\delta(k-1)}{2} \right), \\
    \tau_r(\bar{\delta}(k),\Delta\delta(k)) &= \alpha \frac{d}{2}\left(T(\delta_L(k-1)) - T(\delta_R(k-1))\right) + \nonumber\\
       &\quad \beta \frac{d}{2}\Bigg((a_{j}^{L}-a_{j}^{R})\left(\bar{\delta}(k-1)^2+\frac{\Delta\delta(k-1)^2}{4}\right)\nonumber\\
     &\quad+(a_{j}^{L}+a_{j}^{R})\bar{\delta}(k-1)\Delta\delta(k-1) \nonumber\\
     &\quad + (b_{j}^{L}-b_{j}^{R})\bar{\delta}(k-1) + (b_{j}^{L}+b_{j}^{R})\frac{\Delta\delta(k-1)}{2}\Bigg), %\label{eq_tau_r_dyn}
\end{align*}
where $a_{j}^{i},b_{j}^{i}$ with $i\in\{L,R\},j\in\{f,r\}$ are the forward and reverse parameters to be identified for the propellers. 

%From~\eqref{eq_IG_definition}, 
Applying a time shift as per~\eqref{eq_prop_dyn_model} results in the new equations for the surge, sway, and yaw that yields the following:
\begin{align*}
u(k+1)-u(k) &=\alpha u(k)- \alpha u(k-1)-\alpha h\sigma_u(k-1)\nonumber\\
&\quad+h\bm{M}^{-1}(1,1)\beta \left((a_{j}^{L}+a_{j}^{R})\left(\bar{\delta}(k-1)^2+\frac{\Delta\delta(k-1)^2}{4}\right) \right. \nonumber\\
       &\quad \left.+(a_{j}^{L}-a_{j}^{R})\bar{\delta}(k-1)\Delta\delta(k-1) + (b_{j}^{L}+b_{j}^{R})\bar{\delta}(k-1) \right. \nonumber\\
       &\quad \left. +(b_{j}^{L}-b_{j}^{R})\frac{\Delta\delta(k-1)}{2} \right) + h\sigma_{u}(k), %\label{eq_uk2_uk1_sigma_dyn} 
\end{align*}
\begin{align*}
       v(k+1)-v(k) &=\alpha v(k)- \alpha v(k-1)-\alpha h\sigma_v(k-1)\nonumber\\
&\quad+h\bm{M}^{-1}(2,3)\beta \left((a_{j}^{L}-a_{j}^{R})\left(\bar{\delta}(k-1)^2+\frac{\Delta\delta(k-1)^2}{4}\right) \right. \nonumber\\
       &\quad \left.+(a_{j}^{L}+a_{j}^{R})\bar{\delta}(k-1)\Delta\delta(k-1) + (b_{j}^{L}-b_{j}^{R})\bar{\delta}(k-1) \right. \nonumber\\
       &\quad \left. +(b_{j}^{L}+b_{j}^{R})\frac{\Delta\delta(k-1)}{2} \right) + h\sigma_{v}(k), %\label{eq_vk2_vk1_sigma_dyn}
\end{align*}
\begin{align*}
       r(k+1)-r(k) &=\alpha r(k)- \alpha r(k-1)-\alpha h\sigma_r(k-1)\nonumber\\
&\quad+h\bm{M}^{-1}(3,3)\beta \left((a_{j}^{L}-a_{j}^{R})\left(\bar{\delta}(k-1)^2+\frac{\Delta\delta(k-1)^2}{4}\right) \right. \nonumber\\
       &\quad \left.+(a_{j}^{L}+a_{j}^{R})\bar{\delta}(k-1)\Delta\delta(k-1) + (b_{j}^{L}-b_{j}^{R})\bar{\delta}(k-1) \right. \nonumber\\
       &\quad \left. +(b_{j}^{L}+b_{j}^{R})\frac{\Delta\delta(k-1)}{2} \right) + h\sigma_{r}(k). %\label{eq_rk2_rk1_sigma_dyn} 
\end{align*}

The regressor equations for the surge in matrix form~\eqref{eq:surgedyn}, with the propellers working in the region $(f,f)$, can be derived as follows:
\begin{align*}
    X_{dyn,u} &= \Bigg[\left( \alpha+hR_{u}^{(1)} \right)~~-\alpha hP_{u}^{(1,1)}~~-\alpha hQ_{u}^{(2,3)}~~-\alpha hQ_{u}^{(3,3)}\nonumber \\ &\quad~~- \alpha\left( 1+hR_{u}^{(1)} \right)~~hP_{u}^{(1,1)}~~hQ_{u}^{(2,3)}~~hQ_{u}^{(3,3)}\nonumber \\
        &\quad~~(h-\alpha)\bar{c_{u}}~~2h\beta\bm{M}^{-1}(1,1)a_f~~2h\beta\bm{M}^{-1}(1,1) b_f\Bigg]^T, \\
    A_{dyn,u}(k) &= \Bigg[u(k)~~u(k-1)|u(k-1)|~~v(k-1)r(k-1)~~r(k-1)^2\nonumber \\
        &\quad~~u(k-1)~~u(k)|u(k)|~~v(k)r(k)~~r(k)^2~~1\nonumber \\
        &\quad~~\left(\bar{\delta}(k-1)^2+\frac{\Delta\delta(k-1)^2}{4}\right)~~\bar{\delta}(k-1)\Bigg], \\
    b_{dyn,u}(k) &= u(k+1) - u(k).
\end{align*}

Similarly for sway and yaw, using the dynamic model of the propellers,~\eqref{eq:swaydyn} and~\eqref{eq:yawdyn}, and considering the 3 zones $(f,f)$, $(f,r)$, and $(r,f)$, the following equations represent the regressors for both velocities:
\begin{align*}
       X_{dyn,v} &= \Bigg[\left( \alpha+hR_{v}^{(2)} \right)~~-\alpha hP_{v}^{(2,2)}~~-\alpha hP_{v}^{(2,3)}~~-\alpha hP_{v}^{(3,2)}\nonumber \\ &\quad~~-\alpha hP_{v}^{(3,3)}~~-2\alpha hQ_{v}^{(1,2)}~~-2\alpha hQ_{v}^{(1,3)}~~-\alpha \left(1+hR_{v}^{(2)}\right)\nonumber \\ &\quad~~-\alpha hR_{v}^{(3)}~~ hP_{v}^{(2,2)}~~hP_{v}^{(2,3)}~~hP_{v}^{(3,2)}~~hP_{v}^{(3,3)}\nonumber \\ &\quad~~2hQ_{v}^{(1,2)}~~2hQ_{v}^{(1,3)}~~hR_{v}^{(2)}~~hR_{v}^{(3)}~~(h-\frac{d}{2}\alpha)\bar{c_{v}}\nonumber \\ &\quad~~\frac{d}{2}h\beta\bm{M}^{-1}(2,3)(a_j^L-a_j^R)~~\frac{d}{2}h\beta\bm{M}^{-1}(2,3)(a_j^L+a_j^R)\nonumber \\ &\quad~~\frac{d}{2}h\beta\bm{M}^{-1}(2,3)(b_j^L-b_j^R)~~\frac{d}{2}h\beta\bm{M}^{-1}(2,3)(b_j^L+b_j^R)\Bigg]^T, \\
    A_{dyn,v}(k) &= \Bigg[v(k)~~v(k-1)|v(k-1)|~~v(k-1)|r(k-1)|~~r(k-1)|v(k-1)|\nonumber \\ &\quad~~r(k-1)|r(k-1)|~~u(k-1)v(k-1)~~u(k-1)r(k-1)\nonumber \\ &\quad~~v(k-1)~~r(k-1)~~v(k)|v(k)|~~v(k)|r(k)|~~r(k)|v(k)|\nonumber \\ &\quad~~r(k)|r(k)|~~u(k)v(k)~~u(k)r(k)~~r(k)~~1\nonumber \\ &\quad~~\left(\bar{\delta}(k-1)^2+\frac{\Delta\delta(k-1)^2}{4}\right)~~\bar{\delta}(k-1)\Delta{\delta(k-1)}\nonumber \\ &\quad~~\bar{\delta}(k-1)~~\frac{\Delta{\delta(k-1)}}{2}\Bigg], \\
    b_{dyn,v}(k) &= v(k+1) - v(k).
\end{align*}
\begin{align*}
    X_{dyn,r} &= \Bigg[\left(\alpha+hR_{r}^{(3)}\right)~~-\alpha hP_{r}^{(2,2)}~~-\alpha hP_{r}^{(2,3)}~~-\alpha hP_{r}^{(3,2)}\nonumber \\ &\quad~~-\alpha hP_{r}^{(3,3)}~~-2\alpha hQ_{r}^{(1,2)}~~-2\alpha hQ_{r}^{(1,3)}~~-\alpha hR_{r}^{(2)}\nonumber \\ &\quad~~-\alpha\left(1+hR_{r}^{(3)}\right)~~hP_{r}^{(2,2)}~~hP_{r}^{(2,3)}~~hP_{r}^{(3,2)}~~hP_{r}^{(3,3)}\nonumber \\ &\quad~~2hQ_{r}^{(1,2)}~~2hQ_{r}^{(1,3)}~~hR_{r}^{(2)}~~(h-\frac{d}{2}\alpha)\bar{c_{r}}\nonumber \\ &\quad~~\frac{d}{2}h\beta\bm{M}^{-1}(3,3)(a_j^L-a_j^R)~~\frac{d}{2}h\beta\bm{M}^{-1}(3,3)(a_j^L+a_j^R)\nonumber \\ &\quad~~\frac{d}{2}h\beta\bm{M}^{-1}(3,3)(b_j^L-b_j^R)~~\frac{d}{2}h\beta\bm{M}^{-1}(3,3)(b_j^L+b_j^R)\Bigg]^T, \\
    A_{dyn,r}(k) &= \Bigg[r(k)~~v(k-1)|v(k-1)|~~v(k-1)|r(k-1)|~~r(k-1)|v(k-1)|\nonumber \\ &\quad~~r(k-1)|r(k-1)|~~u(k-1)v(k-1)~~u(k-1)r(k-1)\nonumber \\ &\quad~~v(k-1)~~r(k-1)~~v(k)|v(k)|~~v(k)|r(k)|~~r(k)|v(k)|\nonumber \\ &\quad~~r(k)|r(k)|~~u(k)v(k)~~u(k)r(k)~~v(k)~~1\nonumber \\ &\quad~~\left(\bar{\delta}(k-1)^2+\frac{\Delta\delta(k-1)^2}{4}\right)~~\bar{\delta}(k-1)\Delta{\delta(k-1)}\nonumber \\ &\quad~~\bar{\delta}(k-1)~~\frac{\Delta{\delta(k-1)}}{2}\Bigg], \\
    b_{dyn,r}(k) &= r(k+1) - r(k).
\end{align*}
where $a_{j}^{i},b_{j}^{i}$ with $i\in\{L,R\}$ representing the left and right propellers, and $j\in\{f,r\}$ denote the forward and reverse parameters.

For the sway and yaw velocities in the $(f,f)$ zone, the terms $X_{dyn,v}(19)$ and $X_{dyn,v}(21)$ related to the sway, as well as $X_{dyn,r}(19)$ and $X_{dyn,r}(21)$ related to the yaw, are cancelled.

%% If you have bibdatabase file and want bibtex to generate the
%% bibitems, please use
%%
 % \bibliographystyle{elsarticle-num} 
 % \bibliography{main}

%% else use the following coding to input the bibitems directly in the
%% TeX file.

\end{document}